\documentclass[12pt,a4paper]{article}
\usepackage{inputenc, amssymb, amsmath, amsthm}
\usepackage[english]{babel}
\usepackage{graphicx}
\makeatletter
\thm@style{plain}

\def\th@plain{%
  \thm@headfont{\bfseries}%
  \itshape 
  \thm@notefont{\rm}%
}
\def\thm@indent{\hspace*{\parindent}}
\makeatother

\def\({\left(}
\def\){\right)}

\newcommand{\eps}{\varepsilon}

\newcommand{\be}{\begin{equation}}
\newcommand{\ee}{\end{equation}}

\renewcommand{\leq}{\leqslant}
\renewcommand{\geq}{\geqslant}
\let\epsilon\varepsilon
\let\phi\varphi
\let\le\leqslant
\let\ge\geqslant

\def\arg{\mathop{\ensuremath{\text{\textup{arg}}}}}
\def\sgntwods{\mathop{\ensuremath{\text{\textup{sgn}}_{2\Delta\sigma}}}}

\def\Re{\mathop{\ensuremath{\text{\textup{Re}}}}}
\def\Im{\mathop{\ensuremath{\text{\textup{Im}}}}}

\def\meas{\mathop{\ensuremath{\text{\textup{meas}}}}}

\newtheorem{theo}{Theorem}
\newtheorem{lemma}{Lemma}
\newcommand{\dokvo}{{\it Proof.} }
\newcommand\inte{\int\limits}

\def\bOmega{{\boldsymbol \Omega}}

\newcommand{\beq}{\begin{equation}}
\newcommand{\eeq}{\end{equation}}

{Definition}
{Algorithm}

\newcommand{\ds}{\Delta\sigma}
\newcommand{\K}{\mathcal{K}}
\newcommand{\RR}{\mathcal{R}}
\newcommand{\Ra}{\frac{\alpha-R}{\alpha}}
\newcommand{\aR}{\alpha/(\alpha-R)}

\begin{document}

\centerline{\bf\uppercase{Almost all of the zeros of the Riemann zeta-function}}
\centerline{\bf\uppercase{are on the critical line}\footnote[1]{%
}
}

\bigskip

\medskip

\begin{center}
{\sc Tatyana Preobrazhenskaya}\\[3mm]
{\sc Sergei Preobrazhenskii}\\
\end{center}

\bigskip

\bigskip

\hbox to \textwidth{\hfil\parbox{0.9\textwidth}{%
\small {\sc Abstract.}
This is a reworked version of the paper.

An idea that allows us to circumvent limitations of previous approaches
is not to apply arithmetic-geometric mean inequality and the second moment asymptotics
to the entire segment $[1/2-a/\log T+iT,1/2-a/\log T+i2T]$ but use them on a subset only,
and use the integral of logarithm of the mollified function on the complement.

Ultimately, the result depends on the exponent in the zero-density estimate near the critical line,
which leads to the relation between the magnitude of $\widetilde{V}^{1/2}$ and the measure of the exceptional set
in Theorem 4, Section 2{.}3.
The exponents of Jutila and Conrey are enough for our purposes.

We provide more details on an effective approximation of $1/z$ using the Schwarz--Christoffel mapping.
This is needed in the construction of the mollifier.

One observation on why the approach is feasible is that the functional equation established in the paper
allows one to shift the segment of integration $[1/2-a/\log T+iT,1/2-a/\log T+i2T]$ to $[1/2+A/\log T+iT,1/2+A/\log T+i2T]$.
A result of Selberg allows of a proof that almost all of zeros of our integrand are to the left of the shifted segment.}\hfil}





\section{Introduction}

The Riemann zeta-function $\zeta(s)$ is defined for $\Re s>1$ by
\[
\zeta(s)=\sum_{n=1}^{\infty}n^{-s},
\]
and for other $s$ by the analytic continuation.
It is a meromorphic function in the whole complex plane
with the only singularity $s=1$, which is a simple pole
with residue $1$.

The Euler product links the zeta-function
with prime numbers: for $\Re s>1$
\[
\zeta(s)=\prod_{p\text{\textup{ prime}}}\left(1-p^{-s}\right)^{-1}.
\]

The functional equation for $\zeta(s)$ may be written in the form
\[
\xi(s)=\xi(1-s),
\]
where $\xi(s)$ is an entire function defined by
\[
\xi(s)=H(s)\zeta(s)
\]
with
\[
H(s)=\frac12s(1-s)\pi^{-s/2}\Gamma\left(\frac s2\right).
\]
This implies that $\zeta(s)$ has zeros at
$s=-2$, $-4$, ${\ldots}$
These zeros are called the ``trivial'' zeros.
It is known that $\zeta(s)$ has infinitely many
nontrivial zeros $s=\rho=\beta+i\gamma$, and all of them are in the ``critical strip''
$0<\Re s=\sigma<1$, $-\infty<\Im s=t<\infty$.
The pair of nontrivial zeros with the smallest value of $|\gamma|$
is $\frac12\pm i(14{.}134725\ldots)$.

If $N(T)$ denotes the number of zeros $\rho=\beta+i\gamma$
($\beta$ and $\gamma$ real), for which $0<\gamma\le T$,
then
\[
N(T)=\frac T{2\pi}\log\left(\frac T{2\pi}\right)
-\frac T{2\pi}+\frac78+S(T)+O\left(\frac1T\right),
\]
with
\[
S(T)=\frac1{\pi}\arg\zeta\left(\frac12+iT\right)
\]
and
\[
S(T)=O(\log T).
\]
This is the Riemann--von~Mangoldt formula for $N(T)$.

Let $N_0(T)$ be the number of zeros of $\zeta\left(\frac12+it\right)$ when $0<t\le T$,
each zero counted with multiplicity.
The Riemann hypothesis is the conjecture that $N_0(T)=N(T)$.
Let
\[
\kappa=\liminf\limits_{T\to\infty}\frac{N_0(T)}{N(T)}.
\]

Important results about $N_0(T)$ include:
\begin{itemize}
\item \cite{PreobHardy14}: Hardy proved that $N_0(T)\to\infty$ as $T\to\infty$.
\item \cite{PreobHL21}: Hardy and Littlewood obtained that $N_0(T)\ge AT$ for some $A>0$
and all sufficiently large $T$.
\item \cite{PreobSel42}: Selberg proved that $\kappa\ge A$ for an effectively computable positive constant $A$.
\item \cite{PreobLev74}: Levinson proved that $\kappa\ge0{.}34\ldots$
\item \cite{PreobCon89}: Conrey obtained $\kappa\ge0{.}4088\ldots$
\item \cite{PreobFeng12}: Feng obtained $\kappa\ge0{.}4128\ldots$ (assuming a condition on the lengths of the mollifier)
\end{itemize}

In this article we establish the following
\begin{theo}\label{Preob100percent}
We have
\[
\kappa=1.
\]
\end{theo}

In~\cite{PreobCon83} it is shown
that to estimate the proportion of the critical zeros of the Riemann zeta-function one may use
linear combinations of the $\xi$-function and its derivatives of a fairly general form.
In this paper we select specific linear combinations from them, as per Lemma~\ref{PreobConreyShLemma}.
This lemma asserts that the specific linear combination taken at $s$ is linked to another linear combination of a similar kind,
taken at the point translated by
\[
\Delta\sigma=\frac{\alpha}{\log T}.
\]
It turns out that the possibility of such a translation allows one to improve $\kappa$ substantially,
if one changes some parts of the Levinson--Conrey argument.

We start with a numerical example illustrating the key point of our argument~---
existence of a short Dirichlet polynomial majorant for our specific linear combination
of the Riemann zeta-function and its derivatives.

Let $T=7^{10}$, $\alpha=10$, $R=1$. First,
for $\Re s>1$ we define the function
\[
g_{\alpha,T}(s)\mathrel{:=}-\frac12\sum_{l=1}^{\infty}
\tanh\left(\frac{\alpha}2\left(\frac{\log l}{\log T}-\frac12\right)\right)l^{-s}.
\]

We note that for $s=\frac12-\frac R{\log T}+it$ with $t\in [T/2,T]$,
the analytically continued function
\[
G(s)\mathrel{:=}\frac{\zeta(s)}2+g_{\alpha,T}(s)
\]
can be used as the integrand $G(s)$ in the Levinson--Conrey method
(see Subsection~\ref{PreobOutLevCon}).

\begin{figure}
  \begin{center}
    \includegraphics[width=0.9\textwidth,keepaspectratio]{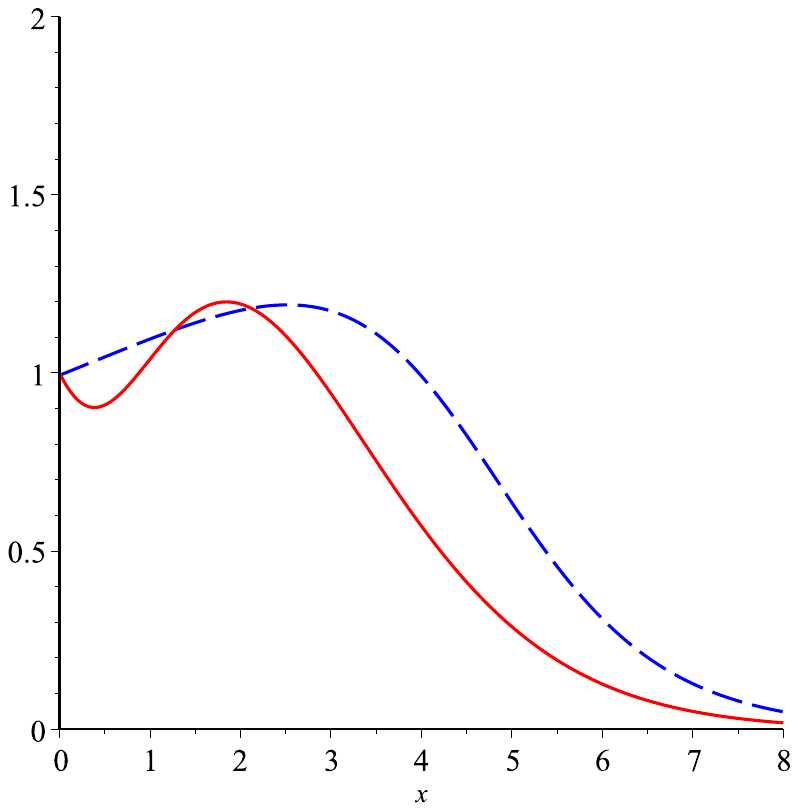}
    \caption{The function
    $a(x)e^{-x}=\frac{e^{\alpha/2}}2\left(1+\tanh\left(\frac 14(2\alpha x/(\alpha-R)-\alpha)\right)\right)e^{-x}$ defining the coefficients
    of the Dirichlet polynomial $\sum_{l\le T}a\left(\frac{(\alpha-R)\log l}{\log T}\right)\exp\left(-\frac{(\alpha-R)\log l}{\log T}\right)l^{-(1/2+it)}$
    (dash), and the approximation $\frac{e^{-x}}{s_6(1)}s_6(x+1)$ (solid)} \label{PreobLaguerre2Fig}
  \end{center}
\end{figure}

Next, for $\Delta\sigma=\frac{\alpha}{\log T}$, using the translation functional equation
\[
2e^{\alpha/2}\left(\frac{\zeta(s+\Delta\sigma)}2-g_{\alpha,T}(s+\Delta\sigma)\right)
=2\left(\frac{\zeta(s)}2+g_{\alpha,T}(s)\right)
\]
(see Sections~\ref{Preobgsection} and~\ref{Preobgsect}) and the approximate functional equation,
we can replace $G(s)$ with the Dirichlet polynomial
\[
\begin{split}
&e^{\alpha/2}\sum_{l\le T}\left(\frac12+\frac12
\tanh\left(\frac{\alpha}2\left(\frac{\log l}{\log T}-\frac12\right)\right)\right)l^{-(s+\Delta\sigma)}\\
=&\frac{e^{\alpha/2}}2\sum_{l\le T}\left(1+\tanh\left(\frac 14\left(\frac{2\alpha\log l}{\log T}-\alpha\right)\right)\right)\exp\left(-\frac{(\alpha-R)\log l}{\log T}\right)l^{-(1/2+it)}
\end{split}
\]
(see Figure~\ref{PreobLaguerre2Fig} (dash)).

\begin{figure}
  \begin{center}
    \includegraphics[width=0.47\textwidth,keepaspectratio]{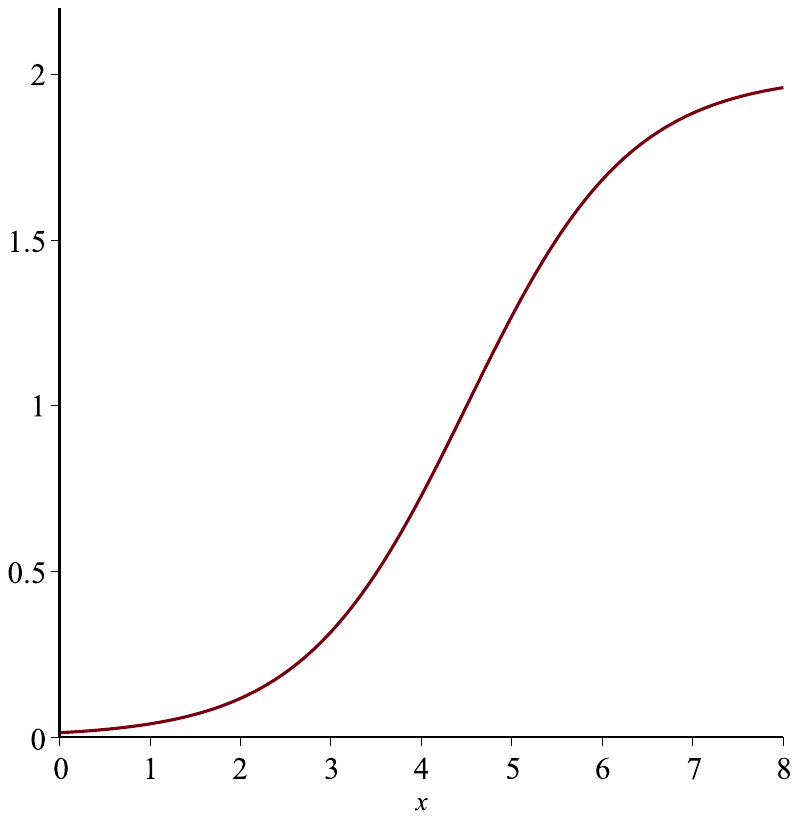}
    \hfill
    \includegraphics[width=0.47\textwidth,keepaspectratio]{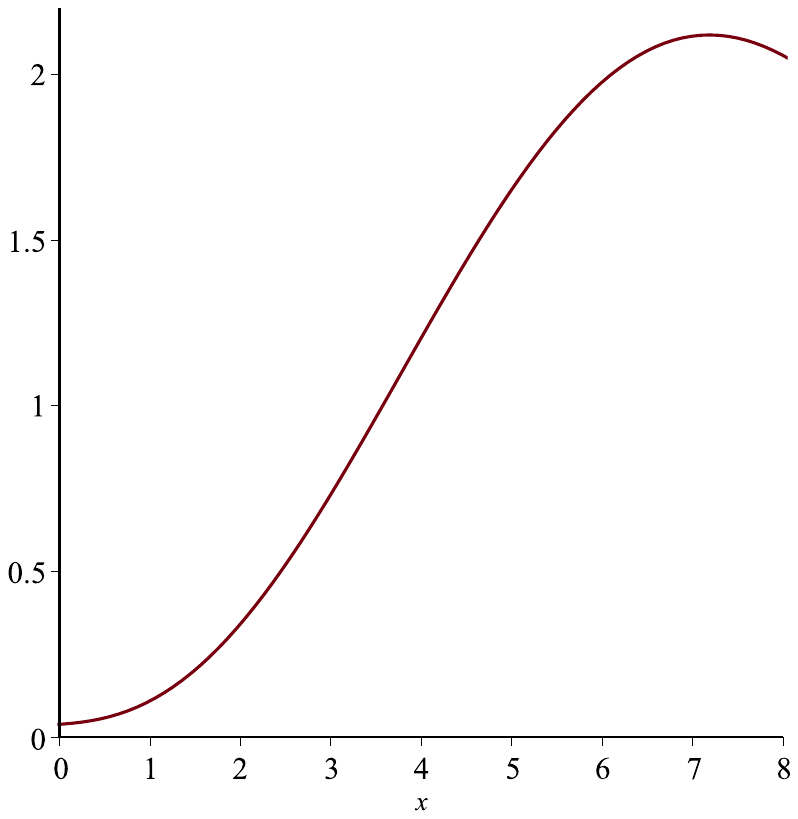}
    \caption{The function $1+\tanh\left(\frac 14(2\alpha x/(\alpha-R)-\alpha)\right)$ and the (shifted) Laguerre polynomial approximation $s_6(x+1)$} \label{PreobLaguerreFig}
  \end{center}
\end{figure}

The next step is to construct the finite Laguerre sum approximation of the function
$1+\tanh\left(\frac 14(2\alpha x/(\alpha-R)-\alpha)\right)$. In this example, we choose the small degree $6$,
and use the shifted Laguerre sum $s_6(x+1)$ (see Figure~\ref{PreobLaguerreFig}), since for this small degree the shifted normalized function
$\frac{e^{-x}}{s_6(1)}s_6(x+1)$ gives better approximation of
\[
\frac{e^{\alpha/2}}2\left(1+\tanh\left(\frac 14(2\alpha x/(\alpha-R)-\alpha)\right)\right)e^{-x}
\]
(see Figure~\ref{PreobLaguerre2Fig}) than the unshifted $\frac{s_6(x)}{s_6(0)}e^{-x}$. For large degrees the shift is unnecessary.

Thus for $\Delta\sigma=\frac{\alpha}{\log T}$, factoring out the Riemann zeta-function $\zeta(s+\Delta\sigma)$, we obtain
the following approximations:
\[
G(s)\approx G^{*}(s+\Delta\sigma)\approx\zeta(s+\Delta\sigma)\left(1+\lambda(s+\Delta\sigma)\right),
\]
where
\[
\lambda(s+\Delta\sigma)
=\frac{c_1(\alpha)}{\log T}\frac{\zeta'}{\zeta}(s+\Delta\sigma)
+\dots+\frac{c_{6}(\alpha)}{(\log T)^{6}}\frac{\zeta^{(6)}}{\zeta}(s+\Delta\sigma)
\]
with the coefficients $c_0(\alpha)=1$, $c_1(\alpha)$, $c_2(\alpha)$, $\dots$, $c_{6}(\alpha)$
defined by the Laguerre sum $\frac1{s_6(1)}s_6(x+1)$ at $x=\frac{(\alpha-R)\log l}{\log T}$.

On the other hand, using the fact that the coefficients of the Dirichlet polynomial
\[
\frac{e^{\alpha/2}}2\sum_{l\le T}\left(1+\tanh\left(\frac 14\left(\frac{2\alpha\log l}{\log T}-\alpha\right)\right)\right)\exp\left(-\frac{(\alpha-R)\log l}{\log T}\right)l^{-(1/2+it)}
\]
are close to $0$ for $\frac{\log l}{\log T}>\frac12$ (see Figure~\ref{PreobLaguerre2Fig})
we get the mean-square estimates for $G(s)$ mollified by the standard mollifier $M(s)$:
(Please note that in the following display we use the oversimplified estimate
$\frac{\gamma(k,(\alpha-R)/2)}{\gamma(k,\alpha-R)}\lesssim\frac1{2^k}$
for the ratio of the lower incomplete gamma functions.)
\[
\begin{split}
&\frac2{T}\inte_{[T/2,T]}\left|G\left(\frac12-\frac R{\log T}+it\right)M\left(\frac12-\frac R{\log T}+it\right)\right|^2dt\\
\approx&\frac2{T}\inte_{[T/2,T]}\left|\sum_{0\le k\le6}\sum_{l\le T^{1/2}}(-1)^kc_k(\alpha)\left(\frac{\log l}{\log T}\right)^k\exp\left(-\frac{(\alpha-R)\log l}{\log T}\right)l^{-(1/2+it)}M\left(\frac12-\frac R{\log T}+it\right)\right|^2dt\\
\lesssim&\frac{14}{T}\sum_{0\le k\le6}\frac1{2^{2k}}\inte_{[T/2,T]}\left|\sum_{l\le T}(-1)^kc_k(\alpha)\left(\frac{\log l}{\log T}\right)^k\exp\left(-\frac{(\alpha-R)\log l}{\log T}\right)l^{-(1/2+it)}M\left(\frac12-\frac R{\log T}+it\right)\right|^2dt\\
\approx&\frac{14}{T}\inte_{[T/2,T]}\sum_{0\le k\le6}\left|\zeta\left(\frac12+\frac{\alpha-R}{\log T}+it\right)M\left(\frac12-\frac R{\log T}+it\right)\right|^2\left|\frac{c_k(\alpha)}{2^k(\log T)^k}\frac{\zeta^{(k)}}{\zeta}\left(\frac12+\frac{\alpha-R}{\log T}+it\right)\right|^2dt,
\end{split}
\]
where
\[
M(s)=\sum_{m\le T^{\theta}}\mu(m)P\left(\frac{\log m}{\log T}\right)m^{-\left(s+\frac R{\log T}\right)},
\]
$\theta>0$ is fixed, $P(x)$ is a real polynomial with $P(0)=1$ and $P(\theta)=0$.

\begin{figure}
  \begin{center}
    \includegraphics[width=0.9\textwidth,keepaspectratio]{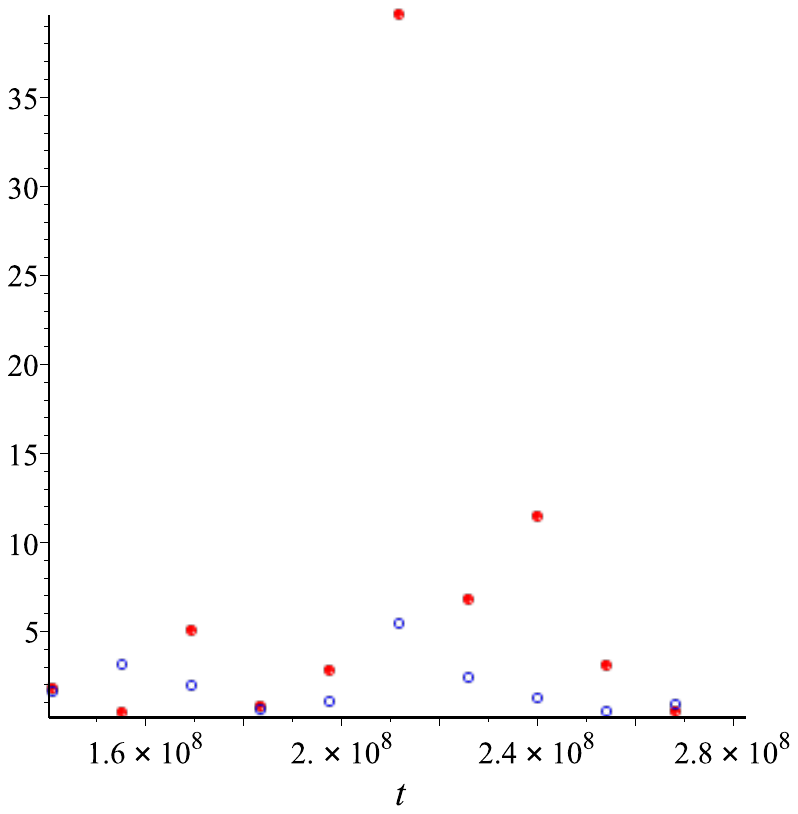}
    \caption{The absolute values of the function $\lambda\left(\frac 12+\frac{\alpha-R}{\log T}+it\right)$ (circles) for $t\in [T/2,T]$ and the majorant $A\left(\frac 12+\frac{\alpha-R}{\log T}+it\right)$ (solid circles) without lower order terms} \label{PreobMajorantFig}
  \end{center}
\end{figure}

The next key step is to provide a short Dirichlet polynomial majorant for
$\lambda(s+\Delta\sigma)$. This is provided by Theorem~\ref{PreobApproximationTh}.
In this example, we put $x=T^{1/10}$, $x_0=T^{1/9}$ and $x_1=T^{1/17}$. This
choice is not fully consistent with our value $\alpha=10$ and $C=1/8$,
but makes the lengths of the Dirichlet polynomials not too large but nontrivial.
In Figure~\ref{PreobMajorantFig} we plot the absolute values of
$\lambda\left(\frac 12+\frac{\alpha-R}{\log T}+it\right)$ (circles) for $t\in [T/2,T]$, and the majorant $A\left(\frac 12+\frac{\alpha-R}{\log T}+it\right)$ (solid circles),
ignoring the terms with $\log s$.

Attaching the mollifying Dirichlet polynomial which is almost the multiplicative inverse
of the majorant $A\left(\frac 12+\frac{\alpha-R}{\log T}+it\right)$, we obtain our functions and the mean-value integrals to be estimated in the Levinson--Conrey framework.

\section{Outline of the method}

\subsection{The Levinson--Conrey Method (A Version of the Principal Inequality)}
\label{PreobOutLevCon}

Given real numbers $\{g_k\}$ let $G(s)$ be the function defined by
\[
G(s)=\frac1{2H(s)}\left(\xi(s)+\sum_{\textup{$k$ odd}}g_k\xi^{(k)}(s)\right)
\]
for $s\in{\mathbb C}$ with $s\ne1$.

The novelty of our specific choice for $G(s)$ given in Subsection~\ref{PreobOutFEq}
is that it obeys the translation functional equation.
It involves higher derivatives of $\zeta$ in an essential way
that pushes almost all of its zeros to the region
\[
\sigma\le1/2-c/\log T
\]
for any fixed $c>0$. A related phenomenon is described in~\cite[Section 3{.}2]{PreobKi11}.
This is why our mollification is so effective.

Next, for $s=\sigma+it$ with $t\asymp T$ and $-1\le\sigma\le2$ we have
\[
G(s)=\sum_{l\le T}Q\left(\frac{\log l}{\log T}+\delta(s)\right)l^{-s}+O\left(T^{-\frac14}\right),
\]
$\delta(s)=\frac{\log(2\pi T/s)}{2\log T}\ll\frac1{\log T}$,
$Q(x)$ is a polynomial such that
\[
Q(x)=\frac12+\frac12\sum\limits_{\substack{k\text{\textup{ odd}}\\ k\le K}}g_k(\log T)^k\left(\frac12-x\right)^k,
\]
with $g_k$ real (or $Q(x)+Q(1-x)=1$).

Consider
\[
F(s)=G(s)M(s){\mathcal L}(A(s))
\]
with
\[
M(s)=\sum_{m\le T^{\theta}}\mu(m)P\left(\frac{\log m}{\log T}\right)m^{-\left(s+\frac R{\log T}\right)},
\]
$\theta>0$ is fixed, $P(x)$ is a real polynomial with $P(0)=1$ and $P(\theta)=0$,
\[
A(s)=\sum_{m\le X}\widetilde{\Lambda}(m)m^{-s}
\]
is a Dirichlet polynomial
with $\widetilde{\Lambda}(1)=0$.
Here $X=T^c$ with $c$ depending on $\alpha$, $\alpha$ is a sufficiently slowly growing function of $T$,
and ${\mathcal L}(w)$ is a polynomial such that
${\mathcal L}(0)=\frac1{Q(0)}$.
Moreover, suppose that the coefficients $c(m)$ of the Dirichlet polynomial $M(s){\mathcal L}(A(s))$
satisfy the bound $c(m)\ll m^{\epsilon}$, where the implied constant can depend on $\alpha$.

The following theorem represents a version of the principal inequality of the Levinson--Conrey method:
\begin{theo}
\label{PreobLevCon}
Let $R$ be fixed, $R>0$, $T$ be a parameter going to infinity,
\[
a=\frac12-\frac R{\log T}.
\]
Let $N_{00}(T,2T)$ be the number of zeros $\rho=\frac12+i\gamma$, $T\le\gamma\le2T$, of $\zeta(s)$
counted without multiplicity, which are not zeros of $G(s)$.
Let $E$ be a subset of $[T,2T]$ which has the measure $\eps_ET$, $0<\eps_E<1$,
and is a union of a finite number of intervals.
Then
\[
N_{00}(T,2T)\ge N(T,2T)\left(1-\frac2{R}(\log I(R)+\eps_E\log I_E(R)+L_E(R))+O\left(\frac1{\log T}\right)\right),
\]
where
\[
I(R)=\frac1{T-\eps_ET}\inte_{[T,2T]\setminus E}|F(a+it)|\,dt,
\]
\[
I_E(R)=\frac1{\eps_ET}\inte_{E}\frac{|F(a+it)|}{|{\mathcal L}(A(a+it))|}\,dt
\]
and
\[
L_E(R)=\frac1{T}\inte_{E}\log|{\mathcal L}(A(a+it))|\,dt.
\]
\end{theo}

Proof of the theorem is given in Section~\ref{princineq}.

\subsection{The Functional Equation}
\label{PreobOutFEq}

\textbf{Definition.}
For $\Re s>1$ define
\[
g_{\alpha,T}(s)\mathrel{:=}-\frac12\sum_{l=1}^{\infty}
\tanh\left(\frac{\alpha}2\left(\frac{\log l}{\log T}-\frac12\right)\right)l^{-s}.
\]

In Sections~\ref{Preobgsection} and~\ref{Preobgsect} we will prove
\begin{enumerate}
\item For $\Delta\sigma=\frac{\alpha}{\log T}$ we have $2e^{\alpha/2}\left(\frac{\zeta(s+\Delta\sigma)}2-g_{\alpha,T}(s+\Delta\sigma)\right)
=2\left(\frac{\zeta(s)}2+g_{\alpha,T}(s)\right)$.
\item For $s=\sigma+it$, $T\le t\le2T$, the function $H(s)g_{\alpha,T}(s)+\text{small perturbation}$
is purely imaginary for $\Re s=\frac12$. The small perturbation term does not affect
the principal inequality of the Levinson--Conrey method as $\alpha$ goes to infinity.
\end{enumerate}

The translation functional equation of item 1 implies
\[
G(s)=\sum_{l\le T}\left(\frac12+q\left(\frac{\log l}{\log T}+\delta(s)\right)\right)l^{-s}+O(T^{-\frac14})=e^{\alpha/2}\sum_{l\le T}\left(\frac12-\widetilde q\left(\frac{\log l}{\log T}+\delta_1(s)\right)\right)l^{-(s+\Delta\sigma)}+{\mathcal R}.
\]
The term $\delta_1(s)=\delta(s+\Delta\sigma)$ and ${\mathcal R}$ comes from
a careful approximation of $-\frac12\tanh\left(\frac{\alpha}2\left(\frac{\log l}{\log T}-\frac12\right)\right)$ by polynomials $q\left(\frac{\log l}{\log T}\right)$ and $\widetilde q\left(\frac{\log l}{\log T}\right)$
of large degrees $K\asymp\alpha$ and $\K\asymp\alpha^3$, respectively~---
see Equations~\eqref{PreobOddPowApprox} and~\eqref{PreobLaguerreApprox} in Section~\ref{Preobgsect} below.

In the right-hand side we denote
\[
G^{*}(s+\Delta\sigma)\mathrel{:=}e^{\alpha/2}\sum_{l\le T}\left(\frac12-\widetilde q\left(\frac{\log l}{\log T}\right)\right)l^{-(s+\Delta\sigma)}.
\]

\begin{theo}
\label{PreobTranslationTh}
In Theorem~\textup{\ref{PreobLevCon}} the function
\[
F(s)=G(s)M(s){\mathcal L}(A(s))
\]
can be replaced by
\[
F^{*}(s)=G^{*}(s+\ds)M(s){\mathcal L}(A(s))
\]
with an acceptable error for $\kappa$, i.e. the error goes to $0$
as $\alpha$ and the degrees $K\asymp\alpha$, $\K\asymp\alpha^3$ of the polynomials $q$ and $\widetilde q$ go to infinity
\textup{(}the conditions on $M(s)$ and ${\mathcal L}(w)$ remain unchanged\textup{)}.
\end{theo}

\textbf{Remark.} For $\Delta\sigma=\frac{\alpha}{\log T}$ we can write
\[
\frac{G^{*}(s+\Delta\sigma)+O(T^{-\frac14})}{\zeta(s+\Delta\sigma)}=c_0(\alpha)+\lambda(s+\Delta\sigma),
\]
where $c_0(\alpha)=(1+o(1))\frac{e^{\alpha/2}}2\left(1-\tanh\left(\frac{\alpha}4\right)\right)$ and
\[
\lambda(s+\Delta\sigma)
=\frac{c_1(\alpha)}{\log T}\frac{\zeta'}{\zeta}(s+\Delta\sigma)
+\dots+\frac{c_{\K}(\alpha)}{(\log T)^{\K}}\frac{\zeta^{(\K)}}{\zeta}(s+\Delta\sigma)
\]
with the coefficients $c_1(\alpha)$, $c_2(\alpha)$, $\dots$, $c_{\K}(\alpha)$
implicitly defined by the polynomial $\widetilde q\left(\frac{\log l}{\log T}\right)$.

\subsection{Theorems of Selberg and Lester (A Generalization)}\label{PreobOutSeLester}

\begin{theo}
\label{PreobApproximationTh}
Let $T\le t\le2T$, $x=T^{1/\alpha^{1/2+4C}}$, $x_0=T^{1/\alpha^{4C}}$, $x_1=T^{1/\alpha^{1/2+8C}}$,
\[
\Delta\sigma=\frac{\alpha}{\log T},
\]
with $\alpha$ going to infinity with $T$ sufficiently slowly,
$R=\epsilon\log\alpha$, $\epsilon>0$ sufficiently small and $C\ge\frac18$ sufficiently large be fixed,
and $\sigma+\ds=\frac12+\frac{\alpha-R}{\log T}$.
Let $\frac{\gamma(k,(\alpha-R)/2)}{\gamma(k,\alpha-R)}$ be the ratio of the lower incomplete gamma functions.
Then there exists a set ${\mathcal M}_{\alpha}\subset[T,2T]$ with $\meas{\mathcal M}_{\alpha}\ge(1-\alpha^{1-C(1-\epsilon)})T$
such that for each $t\in{\mathcal M}_{\alpha}$
we have the following Selberg estimate for $\lambda(\sigma+\ds+it)$, see~\eqref{Preoblex}\textup{:}
\[
\lambda(\sigma+\ds+it)
\ll\left|\sum_{k=1}^{\K}\frac{|c_k(\alpha)|\gamma(k,(\alpha-R)/2)}{\gamma(k,\alpha-R)(\log T)^k}
\sum_{\substack{R_1\ge0,\ldots,R_k\ge0\\R_1+\dots+kR_k=k}}
\frac{k!}{R_1!\dots R_k!}\prod_{j=1}^k\left(\frac{A_j(\sigma+\ds+it)}{j!}\right)^{R_j}\right|,
\]
where
\[
\begin{split}
A_j(\sigma+\ds+it)
=&\frac{j!(\log T)^j}{(\alpha-R-C\log\alpha)^j}
\sum_{n\le x_0^3}\frac{\Lambda_{x_0}(n)}{\log n}
\left(\frac2{n^{\sigma+\ds+(\alpha^{1/2+4C}-\alpha+R)/\log T+it}}\right.\\
&\left.-\frac1{n^{\sigma+\ds+(\alpha^{1/2+8C}-\alpha+R)/\log T+it}}\right)\\
&+\frac{\alpha^{1/2+4C}}{\log T}\frac{j!(\log T)^j}{(\alpha-R-C\log\alpha)^j}
\left(\sum_{n\le x^3}\frac{\Lambda_x(n)}{n^{\sigma+\ds+(\alpha^{1/2+4C}-\alpha+R)/\log T+it}}\right.\\
&\left.+\sum_{n\le x_0^3}\frac{\Lambda_{x_0}(n)}{n^{\sigma+\ds+(\alpha^{1/2+4C}-\alpha+R)/\log T+it}}
+\sum_{n\le x_0^3}\frac{\Lambda_{x_0}(n)}{n^{\sigma+\ds+(\alpha^{1/2+8C}-\alpha+R)/\log T+it}}\right)\\
&+\frac{\alpha^{1/2+8C}}{\log T}\frac{j!(\log T)^j}{(\alpha-R-C\log\alpha)^j}
\sum_{n\le x_1^3}\frac{\Lambda_{x_1}(n)}{n^{\sigma+\ds+(\alpha^{1/2+8C}-\alpha+R)/\log T+it}}\\
&+e^{-\sqrt{\alpha}}\left.\left(\frac{d}{ds_1}\right)^{j}
\inte_{s_1}^{s_1+\frac{\alpha^{1/2+8C}-\alpha+R}{\log T}}\log s\,ds\right|_{s_1=\sigma+\ds+it}.
\end{split}
\]
\end{theo}
We prove Theorem~\ref{PreobApproximationTh} in Section~\ref{selbergsec}.

Next, for all $t\in[T,2T]$ and $c_k(\alpha)$ given in Section~\ref{PreobFEqSec}, ${\widetilde b}_\nu(\alpha)$ in Lemma~\ref{PreobLCoeffsLem}, we denote
\[
  \begin{split}
    &A(\sigma+\Delta\sigma+it)\\
    &\mathrel{:=}\sum_{k=1}^{\K}\frac{|c_k(\alpha)|\gamma(k,(\alpha-R)/2)}{\gamma(k,\alpha-R)(\log T)^k}
    \sum_{\substack{R_1\ge0,\ldots,R_k\ge0\\R_1+\dots+kR_k=k}}
    \frac{k!}{R_1!\dots R_k!}\\
    &\times\widetilde{V}^{1/2}\left|e^{\alpha/2}\frac{(\alpha-R)^k}{(\alpha-R-C\log\alpha)^k}\sum_{\nu=k}^{\K}{\widetilde b}_\nu(\alpha)\binom{\nu}{\nu-k}\right|^{-1}
    \prod_{j=1}^k\left(\frac{A_j(\sigma+\ds+it)}{j!c\sqrt{\log\alpha}}\right)^{R_j}.
  \end{split}
\]
Due to the ratio of the lower incomplete gamma functions
we can truncate the sum by $\K=\alpha\left(\frac12+\epsilon\right)$
and take $\widetilde{V}^{1/2}=\alpha^{C(1/2+2\epsilon)+c}$, $c=o(C)$.
We remark that the bound $\alpha^{1-C(1-\epsilon)}T$ on the measure of the exceptional set
corresponds to the term $-C\log\alpha$ and is a consequence of the Selberg--Jutila zero density estimate.

The behavior of this Dirichlet polynomial $A(\sigma+\Delta\sigma+it)$
is controlled
by rough analogs (see~\eqref{PreobAExtDisk}, \eqref{PreobANarrowStrip}) of the following theorems of Lester.

Let $\psi(T)=(2\sigma-1)\log T$, and for $\psi(T)\geq 1$ define
$V=V(\sigma)=\frac12\sum_{n=2}^{\infty}\frac{\Lambda^2(n)}{n^{2\sigma}}$,\quad
$\Omega=e^{-10}\min\big(V^{3/2},(\psi(T)/\log\psi(T))^{1/2}\big)$.
Suppose that $\psi(T)\to\infty$ with $T$, and $\psi(T)=o(\log T)$.

\textbf{Theorem (Lester's Theorem for a Rectangle).}\textit{
Let $R$ be a rectangle in $\mathbb C$ whose sides are parallel to the coordinate axes.
Then we have
\[
\mathop{\text{\textup{meas}}}\bigg\{t\in(0,T) : \frac{\zeta'}{\zeta}(\sigma+it)V^{-1/2}\in R\bigg\}
=\frac{T}{2\pi}\iint_R\! e^{-(x^2+y^2)/2}\, dx\, dy+O\left(T\frac{(\mathop{\text{\textup{meas}}}(R)+1)}{\Omega}\right).
\]
}
\textbf{Theorem (Lester's Theorem for a Disk).}\textit{
Let $r$ be a real number such that $r\Omega\geq 1$.
Then we have
\[
\mathop{\text{\textup{meas}}}\bigg\{t\in(0,T) : \bigg|\frac{\zeta'}{\zeta}(\sigma+it)\bigg|\leq\sqrt Vr\bigg\}
=\, T(1-e^{-r^2/2})+O\left(T\left(\frac{r^2+r}{\Omega}\right)\right).
\]
}

In Section~\ref{PreobRungeSection} we construct the Dirichlet polynomial ${\mathcal L}(A(s))$
that approximates the function $\frac{\widetilde{\mathcal M}}{c_0+A(s)}$
for almost all values of $t$. Hence in the term
\[
\begin{split}
&\frac2R\log I(R)=\frac2R\log\left(
\frac1{T-\epsilon_ET}\int\limits_{[T,2T]\setminus E}|G^{*}(\sigma+\Delta\sigma+it)M(\sigma+it){\mathcal L}(A(\sigma+\Delta\sigma+it))|\,dt\right)\\
&=\frac2R\\
&\times\log\left(\frac1{T-\epsilon_ET}\int\limits_{[T,2T]\setminus E}|\zeta(\sigma+\Delta\sigma+it)\bigg(c_0+\lambda(\sigma+\Delta\sigma+it)\bigg)M(\sigma+it){\mathcal L}(A(\sigma+\Delta\sigma+it))|\,dt\right)
\end{split}
\]
of the principal inequality of the Levinson--Conrey method, the product
\[
\bigg(c_0+\lambda(\sigma+\Delta\sigma+it)\bigg){\mathcal L}(A(\sigma+\Delta\sigma+it))
\]
is ${}\ll\widetilde{\mathcal M}$ and we will choose $\frac{\log\widetilde{\mathcal M}}R$ to be small.

Now it remains to estimate the integral
\[
\frac1{T-\epsilon_ET}\int\limits_{[T,2T]\setminus E}|\zeta(\sigma+\Delta\sigma+it)M(\sigma+it)|\,dt
\]
for the translated zeta-function $\zeta(\sigma+\Delta\sigma+it)$ and its optimal mollifier $M(\sigma+it)$ (of length $T^{1/2}$, say).
This is done using the mean-square asymptotics. Since $\sigma+\Delta\sigma=\frac12+\frac{\alpha-R}{\log T}$
and we make $\alpha\to\infty$ (slowly) as $T\to\infty$, this integral is close to $1$.

The remaining terms in the principal inequality of the Levinson--Conrey method,
namely, $\epsilon_E\log I_E(R)$ and $L_E(R)$, are proven to give a negligible contribution.

We now proceed to details of the argument.

\section{Translation lemmas}

\begin{lemma}\label{PreobConreyShLemma}
Let $f(s)$ be an analytic function, $s\in{\mathbb C}$, $\ds\in{\mathbb R}$,
$\K\ge1$ be an odd integer.
Then
\be\label{PreobShIdentity}
\begin{split}
&f(s+\ds)=f(s)+\sum_{\substack{k\text{\textup{ odd}}\\ k\le\K}}\left(g_k(\ds)f^{(k)}(s)+g_k(\ds)f^{(k)}(s+\ds)\right)\\
&+\frac{4(-1)^{(\K+1)/2}(\ds)^{\K+1}}{\pi^{\K+2}}\inte_s^{s+\ds}f^{(\K+2)}(w)\left(\sum_{n=1}^{\infty}\frac1{(2n-1)^{\K+2}}\sin\left(\frac{(2n-1)\pi(s+\ds-w)}{\ds}\right)\right)\,dw,
\end{split}
\ee
where
\[
g_k(\ds)=\frac{4(-1)^{(k-1)/2}(\ds)^k}{\pi^{k+1}}\sum_{n=1}^{\infty}\frac1{(2n-1)^{k+1}}.
\]
\end{lemma}

\dokvo We induct on $\K$. First establish induction base $\K=1$. We have
\[
f(s+\ds)=f(s)+\inte_s^{s+\ds}f'(w)\,dw.
\]
Let $\sgntwods(x)$ be the $2\ds$-periodic real-valued function defined by
\[
\sgntwods(x)=\begin{cases}
1&\text{if $x\in(0,\ds)$},\\
0&\text{if $x=-\ds,0,\ds$},\\
-1&\text{if $x\in(-\ds,0)$}.
\end{cases}
\]
Using the Fourier expansion
\be\label{Preobsgnfe}
\sgntwods(x)=\frac4{\pi}\sum_{n=1}^{\infty}\frac1{2n-1}\sin\left(\frac{(2n-1)\pi x}{\ds}\right)
\ee
we obtain
\[
\begin{split}
f(s+\ds)&=f(s)+\frac4{\pi}\inte_s^{s+\ds}f'(w)\left(\sum_{n=1}^{\infty}\frac1{2n-1}\sin\left(\frac{(2n-1)\pi(s+\ds-w)}{\ds}\right)\right)\,dw.
\end{split}
\]
If the range of integration is split over $I_1=[s,s+\epsilon]$, $I_2=[s+\epsilon,s+\ds-\epsilon]$, $I_3=[s+\ds-\epsilon,s+\ds]$
then the integrals over $I_1$ and $I_3$ go to $0$ as $\epsilon\to0$.
The series~\eqref{Preobsgnfe} converges uniformly in $x\in[\epsilon,\ds-\epsilon]$ so by integrating by parts over $I_2$
\[
\begin{split}
&\frac4{\pi}\inte_{s+\epsilon}^{s+\ds-\epsilon}f'(w)\,d\left(\sum_{n=1}^{\infty}\frac{\ds}{(2n-1)^2\pi}\cos\left(\frac{(2n-1)\pi(s+\ds-w)}{\ds}\right)\right)\\
\mathrel{=}&\frac{4\ds}{\pi^2}\sum_{n=1}^{\infty}\frac1{(2n-1)^2}\Bigl(f'(s+\ds-\epsilon)+f'(s+\epsilon)\Bigr)\\
-&\frac{4\ds}{\pi^2}\inte_{s+\epsilon}^{s+\ds-\epsilon}f''(w)\left(\sum_{n=1}^{\infty}\frac1{(2n-1)^2}\cos\left(\frac{(2n-1)\pi(s+\ds-w)}{\ds}\right)\right)\,dw+\delta_1(\epsilon)\\
\mathrel{=}&g_1(\ds)\Bigl(f'(s)+f'(s+\ds)\Bigr)\\
-&\frac{4\ds}{\pi^2}\inte_s^{s+\ds}f''(w)\,d\left(\sum_{n=1}^{\infty}-\frac{\ds}{(2n-1)^3\pi}\sin\left(\frac{(2n-1)\pi(s+\ds-w)}{\ds}\right)\right)+\delta_2(\epsilon)\\
\mathrel{=}&g_1(\ds)\Bigl(f'(s)+f'(s+\ds)\Bigr)\\
-&\frac{4\ds}{\pi^2}\left(-\inte_s^{s+\ds}-\frac{\ds}{\pi}f'''(w)\left(\sum_{n=1}^{\infty}\frac1{(2n-1)^3}\sin\left(\frac{(2n-1)\pi(s+\ds-w)}{\ds}\right)\right)dw\right)+\delta_3(\epsilon),
\end{split}
\]
and $\delta_1(\epsilon),\delta_2(\epsilon),\delta_3(\epsilon)\to0$ as $\epsilon\to0$.
This proves the induction base. The induction step is proven by integrating by parts in~\eqref{PreobShIdentity}
as above,
with the uniform convergence of the series in the integrand when $\K\ge1$.

\noindent{\bf Remark.} We have
\[
g_1(\ds)=\frac{4\ds}{\pi^2}\sum_{n=1}^{\infty}\frac1{(2n-1)^2}=\frac{4\ds}{\pi^2}\zeta(2)\left(1-\frac1{2^2}\right)=\frac{\ds}2,
\]
and in general for $k$ odd
\be\label{PreobBernoulliEq}
g_k(\ds)=\frac{4(-1)^{(k-1)/2}(\ds)^k}{\pi^{k+1}}\zeta(k+1)\left(1-\frac1{2^{k+1}}\right)
=-\left(\frac{\ds}2\right)^k\frac{2^{k+1}-4^{k+1}}{(k+1)!}B_{k+1},
\ee
where $B_{k+1}$ is the Bernoulli number.

The series
\[
\sum_{\substack{k\ge1\\ k\text{\textup{ odd}}}}\left(g_k(\ds)f^{(k)}(s)+g_k(\ds)f^{(k)}(s+\ds)\right)
\]
obtained by successive integrations by parts in~\eqref{PreobShIdentity}
may be divergent. However, we have the following Lemma~\ref{PreobBernoulliLemma}.
Note that if in the series
\[
\sum_{\substack{k\ge1\\ k\text{\textup{ odd}}}}\Bigl(-g_k(\ds)\Bigr)(\log T)^k\left(\frac12-x\right)^k
\]
considered in Lemma~\ref{PreobBernoulliLemma}, we substitute $x=\frac{\log l}{\log T}$,
multiply over by $l^{-(s+\ds)}$ and sum over $l$ from $1$ to $T$,
then we get an approximation for
\[
H(s+\ds)^{-1}\sum_{\substack{k\ge1\\ k\text{\textup{ odd}}}}\Bigl(-g_k(\ds)\Bigr){\xi}^{(k)}(s+\ds)
\]
(see~\cite[Chapter 18, (18{.}9)]{PreobIw14}).
\begin{lemma}\label{PreobBernoulliLemma}
Suppose that $0<\epsilon<2\pi$, $|\alpha|\le2\pi-\epsilon$
and
\[
|\ds|=\frac{|\alpha|}{\log T}\le\frac{2\pi-\epsilon}{\log T}.
\]
Then the series
\[
\sum_{\substack{k\ge1\\ k\text{\textup{ odd}}}}\Bigl(-g_k(\ds)\Bigr)(\log T)^k\left(\frac12-x\right)^k
\]
converges on $x\in[0,1]$ and
\[
\sum_{\substack{k\ge1\\ k\text{\textup{ odd}}}}\Bigl(-g_k(\ds)\Bigr)(\log T)^k\left(\frac12-x\right)^k
=-\tanh\left(\frac{\alpha}2\left(\frac12-x\right)\right).
\]
\end{lemma}

\dokvo By the definition of the Bernoulli numbers,
\[
\frac{z}{e^z-1}=\sum_{m=0}^{\infty}\frac{B_m}{m!}z^m,
\]
with $B_0=1$, $B_1=-\frac12$ and $B_3=B_5=B_7=\dots=0$,
the radius of convergence of the series being $2\pi$.

Since for $|z|\le\pi/2-\epsilon$
\[
\tanh(z)=\sum_{n\ge1}\frac{2^{2n}\left(2^{2n}-1\right)B_{2n}}{(2n)!}z^{2n-1},
\]
then
\[
-\tanh\left(\frac{\alpha}2\left(\frac12-x\right)\right)
=\sum_{\substack{k\ge1\\ k\text{\textup{ odd}}}}\frac{(2^{k+1}-4^{k+1})B_{k+1}}{(k+1)!}\frac{\alpha^k}{2^k}
\left(\frac12-x\right)^k
\]
and the lemma follows from~\eqref{PreobBernoulliEq}.

The following lemma is an easy consequence of Stirling's formula.
\begin{lemma}\label{PreobHsplusdsLemma}
Suppose that $\alpha$ is real
and
\[
\ds=\frac{\alpha}{\log T}.
\]
Then in the rectangle
\[
s=\sigma+it,\quad\frac13\le\sigma\le A,\quad T\le t\le2T
\]
with $A\ge3$ and $T\ge2A$
we have
\[
H(s+\ds)=\left(e^{\alpha/2}+O_{\alpha}\left(\frac1{\log T}\right)\right)H(s).
\]
\end{lemma}
Thus Lemmas~\ref{PreobConreyShLemma} (with Remark), \ref{PreobBernoulliLemma}, \ref{PreobHsplusdsLemma}
allow one to replace
\[
e^{\alpha/2}\sum_{l=1}^{T}\left(1+\tanh\left(\frac{\alpha}2\left(\frac{\log l}{\log T}-\frac12\right)\right)\right)l^{-(s+\ds)}
\]
by
\[
\sum_{l=1}^{T}\left(1-\tanh\left(\frac{\alpha}2\left(\frac{\log l}{\log T}-\frac12\right)\right)\right)l^{-s}
\]
plus some error,
i.e. to link the value of the function suitable for the Levinson--Conrey method at $s$
to the value of the similarly looking function at $s+\ds$,
yet subject to $|\alpha|<2\pi$.
In~Section~\ref{Preobgsection} we shall get rid of this constraint.

\section{A version of the principal inequality of the Levinson--Conrey method}\label{princineq}

For a detailed exposition of the Levinson--Conrey method, see~\cite{PreobCon89}, \cite{PreobIw14}.

Suppose that in the rectangle
\[
s=\sigma+it,\quad\frac13\le\sigma\le A,\quad T\le t\le2T,
\]
with $A\ge3$ and $T\ge2A$, function $G(s)$ is of the form
\[
G(s)=\sum_{l\le T}Q\left(\frac{\log l}{\log T}+\delta(s)\right)l^{-s}+O\left(T^{-\frac14}\right),
\]
where
\[
\delta(s)=\frac{\log(2\pi T/s)}{2\log T}\ll\frac1{\log T},
\]
and $Q(x)$ is a polynomial such that
\[
Q(x)=\frac12+\frac12\sum_{\substack{k\text{\textup{ odd}}\\ k\le\K}}g_k(\log T)^k\left(\frac12-x\right)^k,
\]
where $g_k$ are real (the equation is equivalent to $Q(x)+Q(1-x)=1$).

Next, consider
\[
F(s)=G(s)M(s){\mathcal L}(A(s))
\]
with the mollifiers
\[
M(s)=\sum_{m\le T^{\theta}}\mu(m)P\left(\frac{\log m}{\log T}\right)m^{-\left(s+\frac R{\log T}\right)},
\]
where $P(x)$ is a real polynomial with
\[
P(0)=1\qquad\text{and}\qquad P(\theta)=0,
\]
$A(s)$ is a Dirichlet polynomial
\[
A(s)=\sum_{m\le X}\widetilde{\Lambda}(m)m^{-s}
\]
with $\widetilde{\Lambda}(1)=0$
(here $X=T^c$ with $c$ depending on $\alpha$),
and ${\mathcal L}(w)$ is a polynomial such that
${\mathcal L}(0)=\frac1{Q(0)}$.
Moreover, suppose that the coefficients $c(m)$ of the Dirichlet polynomial $M(s){\mathcal L}(A(s))$
satisfy the bound $c(m)\ll m^{\epsilon}$, where the implied constant can depend on $\alpha$.

We now prove Theorem~\ref{PreobLevCon}, Subsection~\ref{PreobOutLevCon}
which represents a version of the principal inequality of the Levinson--Conrey method.

\textbf{Proof of Theorem~\ref{PreobLevCon}.} The inequality
\[
N_{00}(T,2T)\ge N(T,2T)\left(1-\frac2{R}l(R)+O\left(\frac1{\log T}\right)\right)
\]
where $l(R)=\frac1{T}\int_T^{2T}\log|F(a+it)|\,dt$
is provided in~\cite[Chapter 22, (22{.}4) and (22{.}5)]{PreobIw14}.

Our function $F(s)$ has additional mollifier ${\mathcal L}(A(s))$.
The difference between our mollifier and the one in the book is that our constants can depend on $\alpha$.
But in our argument we can suppose that $\alpha$ is fixed.
The larger $\alpha$ we take, the closer we approach $\varkappa=1$ in the end.
Eventually we can take $\alpha$ to be a sufficiently slowly growing function of $T$.

Next we write
\[
\begin{split}
&\frac1{T}\int_T^{2T}\log|F(a+it)|\,dt
\le\frac1{T-\epsilon_ET}\int\limits_{[T,2T]\setminus E}\log|F(a+it)|\,dt\\
&+\frac{\epsilon_E}{\epsilon_ET}\int\limits_{E}\log\frac{|F(a+it)|}{|{\mathcal L}(A(a+it))|}\,dt
+\frac1{T}\int\limits_{E}\log|{\mathcal L}(A(a+it))|\,dt.
\end{split}
\]
The inequalities
\[
\frac1{T-\epsilon_ET}\int\limits_{[T,2T]\setminus E}\log|F(a+it)|\,dt
\le\log\left(\frac1{T-\epsilon_ET}\int\limits_{[T,2T]\setminus E}|F(a+it)|\,dt\right)
\]
and
\[
\frac1{\epsilon_ET}\int\limits_{E}\log\frac{|F(a+it)|}{|{\mathcal L}(A(a+it))|}\,dt
\le\log\left(\frac1{\epsilon_ET}\int\limits_{E}\frac{|F(a+it)|}{|{\mathcal L}(A(a+it))|}\,dt\right)
\]
follow by considering the integral sums and using arithmetic--geometric mean inequality.

\section{Function $g_{\alpha,T}(s)$}\label{Preobgsection}
In the subsequent arguments we shall get rid of the limitation
$|\alpha|<2\pi$ in Lemma~\ref{PreobBernoulliLemma} by showing
that for $\alpha$ arbitrarily large
the analytic function given for $\Re s>1$ by
\[
g_{\alpha,T}(s)=-\frac12\sum_{l=1}^{\infty}\tanh\left(\frac{\alpha}2\left(\frac{\log l}{\log T}-\frac12\right)\right)l^{-s}
\]
obeys two types of symmetries:
\begin{enumerate}
\item For $\ds=\frac{\alpha}{\log T}$
\be\label{PreobS1}
2e^{\alpha/2}\left(\frac{\zeta(s+\ds)}2-g_{\alpha,T}(s+\ds)\right)=2\left(\frac{\zeta(s)}2+g_{\alpha,T}(s)\right).
\ee
\item For $s=\sigma+it$, $T\le t\le2T$, the function
\be\label{PreobS2}
H(s)g_{\alpha,T}(s)
\ee
is approximated by a sum of the odd derivatives of the $\xi$ function with real coefficients
which are purely imaginary for $\Re s=\frac12$.
\end{enumerate}

To prove~\eqref{PreobS1} we note that
\[
e^{-u}\left(1+\tanh\left(\frac u2\right)\right)=1-\tanh\left(\frac u2\right).
\]
Now take
\[
u=\alpha\left(\frac{\log l}{\log T}-\frac12\right)
\]
multiply the first formula by $l^{-s}$ with $\Re s>1$ and sum over $l\ge1$.

We shall prove~\eqref{PreobS2} in Lemma~\ref{PreobApproxHgLem}.

In the following section, we shall describe properties of $g_{\alpha,T}(s)$ in detail.

\section{Properties of $g_{\alpha,T}(s)$}\label{Preobgsect}

\begin{lemma}[Analytic continuation of $g_{\alpha,T}(s)$]
\label{PreobAnContgLem}
For $s=\sigma+it$ with $\sigma>0$ and $0<t_0\le|t|\le2T$, where $t_0$ is fixed and $T\ge1$,
and for integer $N\ge T$ we have
\[
\begin{split}
g_{\alpha,T}(s)=-\frac12&\left(\sum_{n=1}^N\tanh\left(\frac{\alpha}2\left(\frac{\log n}{\log T}-\frac12\right)\right)n^{-s}\right.\\
&\left.{}-\frac{2T^{1-s}\log T}{\alpha(e^{\alpha/2}+1)(1-(1-s)(\log T)/\alpha)}F(1,1;2-(1-s)(\log T)/\alpha;(e^{\alpha/2}+1)^{-1})\right.\\
&\left.{}+\frac{T^{1-s}}{s-1}-\inte_T^{N+1/2}\tanh\left(\frac{\alpha}2\left(\frac{\log u}{\log T}-\frac12\right)\right)u^{-s}du\right.\\
&\left.{}+\frac{\alpha}{2\log T}\inte_{N+1/2}^{+\infty}\psi(u)\cosh^{-2}\left(\frac{\alpha}2\left(\frac{\log u}{\log T}-\frac12\right)\right)u^{-s-1}du\right.\\
&\left.{}-s\inte_{N+1/2}^{+\infty}\psi(u)\tanh\left(\frac{\alpha}2\left(\frac{\log u}{\log T}-\frac12\right)\right)u^{-s-1}du\right),
\end{split}
\]
where $F(a,b;c;z)$ is the hypergeometric function, and $\psi(x)=x-[x]-\frac12$.
\end{lemma}

\dokvo By the exact summation formula we have
\[
\begin{split}
&\sum_{N+1/2<n\le M+1/2}\tanh\left(\frac{\alpha}2\left(\frac{\log n}{\log T}-\frac12\right)\right)n^{-s}
=\inte_{N+1/2}^{M+1/2}\tanh\left(\frac{\alpha}2\left(\frac{\log u}{\log T}-\frac12\right)\right)u^{-s}du\\
{}+&\frac{\alpha}{2\log T}\inte_{N+1/2}^{M+1/2}\psi(u)\cosh^{-2}\left(\frac{\alpha}2\left(\frac{\log u}{\log T}-\frac12\right)\right)u^{-s-1}du\\
{}-&s\inte_{N+1/2}^{M+1/2}\psi(u)\tanh\left(\frac{\alpha}2\left(\frac{\log u}{\log T}-\frac12\right)\right)u^{-s-1}du.
\end{split}
\]
The first integral is convergent for $\sigma>1$ as $M\to+\infty$,
whereas the latter two integrals with the $\psi$ function converge absolutely for $\sigma>0$.
Denote them by $\Psi_1$ and $\Psi_2$.
Now for $\sigma>1$ we have the formula
\[
\begin{split}
g_{\alpha,T}(s)=-\frac12&\left(\sum_{n=1}^N\tanh\left(\frac{\alpha}2\left(\frac{\log n}{\log T}-\frac12\right)\right)n^{-s}\right.\\
&\left.{}+\inte_T^{+\infty}\tanh\left(\frac{\alpha}2\left(\frac{\log u}{\log T}-\frac12\right)\right)u^{-s}du
-\inte_T^{N+1/2}\tanh\left(\frac{\alpha}2\left(\frac{\log u}{\log T}-\frac12\right)\right)u^{-s}du\right.\\
&\left.{}+\frac{\alpha}{2\log T}\Psi_1-s\Psi_2\right),
\end{split}
\]
in which we consider
\[
\inte_T^{+\infty}\tanh\left(\frac{\alpha}2\left(\frac{\log u}{\log T}-\frac12\right)\right)u^{-s}du.
\]
We write the integrand as
\[
\tanh\left(\frac{\alpha}2\left(\frac{\log u}{\log T}-\frac12\right)\right)u^{-s}
=-\frac{2u^{-s}}{e^{-\alpha/2}u^{\alpha/\log T}+1}+u^{-s}.
\]
Integrating the latter term we get $\frac{T^{1-s}}{s-1}$,
while the former term gives
\[
\begin{split}
\inte_T^{+\infty}\frac{-2T^{-\frac s2+\frac12}}{\left(\frac u{\sqrt T}\right)^{\frac{\alpha}{\log T}}+1}
\left(\frac u{\sqrt T}\right)^{-s}d\left(\frac u{\sqrt T}\right)
&=-2T^{\frac{1-s}2}\inte_{\sqrt T}^{+\infty}\frac{x^{-s}dx}{1+x^{\alpha/\log T}}\\
&=\frac{-2T^{\frac{1-s}2}\log T}{\alpha}\inte_{e^{\alpha/2}}^{+\infty}\frac{v^{(1-s)(\log T)/\alpha-1}}{1+v}\,dv.
\end{split}
\]
Making the change of variables
\[
\begin{split}
w&=\frac1{v+1},\\
v&=\frac1w-1,\\
dv&=-\frac1{w^2}\,dw
\end{split}
\]
we get the integral
\[
\frac{-2T^{\frac{1-s}2}\log T}{\alpha}\inte_0^{(e^{\alpha/2}+1)^{-1}}w^{1-(1-s)(\log T)/\alpha-1}(1-w)^{(1-s)(\log T)/\alpha-1}dw
\]
that can be written as the incomplete beta function
\[
\frac{-2T^{\frac{1-s}2}\log T}{\alpha}B_{(e^{\alpha/2}+1)^{-1}}(1-(1-s)(\log T)/\alpha,(1-s)(\log T)/\alpha)
\]
which in turn can be expressed in terms of the hypergeometric function
\[
\begin{split}
&\frac{-2T^{\frac{1-s}2}\log T}{\alpha}(e^{\alpha/2}+1)^{(1-s)(\log T)/\alpha}\\
&\times F\left(1-(1-s)(\log T)/\alpha,1-(1-s)(\log T)/\alpha;2-(1-s)(\log T)/\alpha;(e^{\alpha/2}+1)^{-1}\right).
\end{split}
\]
Using the known linear transformation formula
\[
F(a,b;c;z)=(1-z)^{c-a-b}F(c-a,c-b;c;z)
\]
we get the term
\[
-\frac{2T^{1-s}\log T}{\alpha(e^{\alpha/2}+1)(1-(1-s)(\log T)/\alpha)}
F(1,1;2-(1-s)(\log T)/\alpha;(e^{\alpha/2}+1)^{-1})
\]
of the analytic continuation formula,
where the function
\[
F(1,1;2-(1-s)(\log T)/\alpha;(e^{\alpha/2}+1)^{-1})
\]
is analytic and bounded in $s$ for $|t|\ge t_0>0$ by the series representation.

\begin{lemma}[Approximate equation for $g_{\alpha,T}(s)$]
\label{PreobApproxEqgLem}
For $s=\sigma+it$ with $\sigma\ge\sigma_0>0$ and $0<t_0\le|t|\le2T$, where $\sigma_0$, $t_0$ are fixed and $T\ge1$,
we have
\[
\begin{split}
g_{\alpha,T}(s)=-\frac12&\left(\sum_{n=1}^T\tanh\left(\frac{\alpha}2\left(\frac{\log n}{\log T}-\frac12\right)\right)n^{-s}\right.\\
&\left.{}-\frac{2T^{1-s}\log T}{\alpha(e^{\alpha/2}+1)(1-(1-s)(\log T)/\alpha)}F(1,1;2-(1-s)(\log T)/\alpha;(e^{\alpha/2}+1)^{-1})\right.\\
&\left.{}+\frac{T^{1-s}}{s-1}+O\left(T^{-\sigma}\right)\right),
\end{split}
\]
where the constant in the $O$-term is absolute.
\end{lemma}

\dokvo In the analytic continuation formula of Lemma~\ref{PreobAnContgLem}
we use the standard uniform approximation
\[
\sum_{T<n\le N+1/2}\tanh\left(\frac{\alpha}2\left(\frac{\log n}{\log T}-\frac12\right)\right)n^{-s}
=\inte_T^{N+1/2}\tanh\left(\frac{\alpha}2\left(\frac{\log u}{\log T}-\frac12\right)\right)u^{-s}du+O\left(T^{-\sigma}\right)
\]
and make $N\to\infty$. $\qed$

We now obtain approximations to $g_{\alpha,T}(s)$
that we need in the context of Conrey's construction~\cite[Chapter 18]{PreobIw14}.
First, we approximate it by using the Fourier expansion
\[
\tanh\left(\frac{\alpha x}2\right)=\sum_{k=1}^Kb_k(\alpha)\sin(kx)+\hat{R}_{K,\alpha}(x)
\]
and the Taylor expansion
\[
\sin(kx)=\sum_{m=1}^M(-1)^{m-1}\frac{(kx)^{2m-1}}{(2m-1)!}+R_{k,M}(x),
\]
where
\[
\begin{split}
\hat{R}_{K,\alpha}(x)&=-2\inte_0^{\pi}\phi_{\alpha,x}(y)D_K(y)\,dy,\\
\phi_{\alpha,x}(y)&=\frac{\tanh\left(\frac{\alpha(x+y)}2\right)+\tanh\left(\frac{\alpha(x-y)}2\right)-2\tanh\left(\frac{\alpha x}2\right)}2,\\
R_{k,M}(x)&=\frac{(-1)^M(kx)^{2M+1}}{(2M)!}\inte_0^1(1-u)^{2M}\cos(kxu)\,du,
\end{split}
\]
and $D_K(y)$ is the Dirichlet kernel.
Explicitly, the coefficients $b_k(\alpha)$ are
\[
\begin{split}
b_k(\alpha)&=-\frac4{\pi}\inte_0^{\pi}\frac{e^{ikx}-e^{-ikx}}{(e^{\alpha x}+1)2i}\,dx\\
&=-\frac4{\pi\alpha 2i}\left(\inte_1^{e^{\alpha\pi}}\frac{v^{ik/\alpha-1}}{v+1}\,dv-\inte_1^{e^{\alpha\pi}}\frac{v^{-ik/\alpha-1}}{v+1}\,dv\right)\\
&=-\frac4{\pi\alpha}\Im\left(B_{1/2}\left(1-i\frac k{\alpha},i\frac k{\alpha}\right)-B_{(e^{\alpha\pi}+1)^{-1}}\left(1-i\frac k{\alpha},i\frac k{\alpha}\right)\right),
\end{split}
\]
where $B_x(a,b)$ is the incomplete beta function.

So we have
\be
\label{PreobOddPowApprox}
\begin{split}
&\sum_{l=1}^T\tanh\left(\frac{\alpha}2\left(\frac{\log l}{\log T}-\frac12\right)\right)l^{-s}\\
&=\sum_{k=1}^Kb_k(\alpha)\sum_{m=1}^M(-1)^{m-1}\frac{k^{2m-1}}{(2m-1)!}
\sum_{l\le T}\left(\frac{\log l}{\log T}-\frac12\right)^{2m-1}l^{-s}\\
&+\sum_{l\le T}{\mathcal R}_{K,M,\alpha}\left(\frac{\log l}{\log T}-\frac12\right)l^{-s},
\end{split}
\ee
where
\[
{\mathcal R}_{K,M,\alpha}(x)=\hat{R}_{K,\alpha}(x)+\sum_{k=1}^Kb_k(\alpha)R_{k,M}(x).
\]
We multiply the polynomial
\[
\sum_{k=1}^Kb_k(\alpha)\sum_{m=1}^M(-1)^{m-1}\frac{k^{2m-1}}{(2m-1)!}
\left(\frac{\log l}{\log T}-\frac12\right)^{2m-1}
\]
appearing in the right-hand side of~\eqref{PreobOddPowApprox} by $-\frac12$ and denote it by
$q\left(\frac{\log l}{\log T}\right)$.

\begin{lemma}[Approximation to $H(s)g_{\alpha,T}(s)$ by a sum of the odd derivatives of the $\xi$ function]
\label{PreobApproxHgLem}
For $s=\sigma+it$ in the rectangle $\frac13\le\sigma\le A$, $T\le t\le2T$, with $A\ge3$ and $T\ge2A$,
we have
\[
2H(s)\left(\sum_{l\le T}\left(q\left(\frac{\log l}{\log T}+\delta(s)\right)
+{\mathcal R}_{K,M,\alpha}\left(\frac{\log l}{\log T}-\frac12\right)\right)l^{-s}+O(T^{-\frac14})\right)
=\sum_{m=1}^Mg_{2m-1}\xi^{(2m-1)}(s),
\]
where
\[
\delta(s)=\frac{\log(2\pi T/s)}{2\log T}\ll\frac1{\log T},
\]
\[
{\mathcal R}_{K,M,\alpha}\ll e^{-\alpha}\quad\text{for some}\quad K\asymp\alpha
\]
and $g_{2m-1}$ are real numbers.
\end{lemma}

\dokvo See~\cite[Chapter 18]{PreobIw14}.

Now in Section~\ref{princineq} we can substitute
\[
Q(x)=\frac12+q(x)
\]
for $Q(x)$ (with $2M-1$ in place of $K$),
and
$\alpha$ can be arbitrarily large.

For the translated function $G^{*}(s+\ds)$ we shall use Laguerre polynomials expansion
of $-\frac12\tanh\left(\frac{\alpha}2\left(y-\frac12\right)\right)$ to construct
$\widetilde q\left(\frac{\log l}{\log T}\right)$, since in this case,
for
\[
\widetilde Q(y)=\frac{1+\tanh\left(\frac{\alpha}2\left(y-\frac12\right)\right)}{1-\tanh\frac{\alpha}4}
\]
the integral $\widetilde A$ corresponding to $A$
in~\eqref{PreobMeanSqInt} will be
\[
\widetilde A=\frac1{(1-\tanh(\alpha/4))^2}\inte_0^1\left(\left(1+\tanh\left(\frac{\alpha}2\left(y-\frac12\right)\right)\right)e^{-(\alpha-R)y}\right)^2dy.
\]
We write the finite Laguerre polynomial expansion of $1+\tanh\left(\frac{\frac{2\alpha}{\alpha-R}x-\alpha}4\right)$
\[
s_{\K-1}(\alpha,x)=\widetilde b_0(\alpha)L_0(x)+\widetilde b_1(\alpha)L_1(x)+\dots+\widetilde b_{\K-1}(\alpha)L_{\K-1}(x).
\]
However,
\[
s_{\K-1}(\alpha,0)
\]
may not be equal to $1-\tanh(\alpha/4)$. We need this condition on $\widetilde q(0)$ in the application of Littlewood's lemma.
But from the representation
\[
s_{\K-1}(\alpha,x)=\inte_0^{+\infty}\left(1+\tanh\left(\frac{\frac{2\alpha}{\alpha-R}t-\alpha}4\right)\right)\widetilde K_{\K-1}(t,x)e^{-t}dt
\]
with the kernel $\widetilde K_{\K-1}(t,x)$ and from~\eqref{PreobRodrig}--\eqref{PreobTanhDer}
it follows in a standard way using the method of steepest descent with the saddle point
$$
\frac{\sqrt{8 \alpha  (\alpha -2 i \pi ) k (R-\alpha )+(\alpha  (\alpha -R-2)+2 i \pi  (R-\alpha ))^2}+2 i \pi  (R-\alpha )+\alpha  (\alpha -R-2)}{4 (2 \pi +i \alpha ) \alpha }
$$
that for a certain $\K\gg\alpha^3$
\[
\left|s_{\K-1}(\alpha,x)
-\left(1+\tanh\left(\frac{\frac{2\alpha}{\alpha-R}x-\alpha}4\right)\right)\right|\le e^{-10\alpha}
\]
uniformly in $x\in[0,\alpha/2]$.
Reverting to the original variable, we denote
\beq
\label{PreobLaguerreApprox}
\widetilde q\left(\frac{\log l}{\log T}\right)=-\frac12\left(s_{\K-1}\left(\alpha,(\alpha-R)\frac{\log l}{\log T}\right)-1\right).
\eeq
Here we record the known formulas:
the Laguerre kernel~\cite[(5{.}1{.}11)]{PreobSz75}
\be
\label{PreobKer}
\widetilde K_{\K-1}(x,u)=\K\frac{L_{\K-1}(x)L_\K(u)-L_\K(x)L_{\K-1}(u)}{x-u},
\ee
the Rodrigues formula~\cite[(5{.}1{.}5)]{PreobSz75}
\be
\label{PreobRodrig}
e^{-u}L_k(u)=\frac1{k!}\frac{d^k}{du^k}\left(e^{-u}u^k\right),
\ee
and the integral representation
\be
\label{PreobTanhDer}
\begin{split}
&\frac14\frac{d^k}{du^k}\left(\tanh\left(\frac14\left(\frac{2\alpha}{\alpha-R}u-\alpha\right)\right)
\right)\\
=&\inte_{0}^{+\infty}\frac{\frac{\partial^k}{\partial u^k}
\left(\sin\left(\left(\frac{2\alpha}{\alpha-R}u-\alpha\right)v\right)
\right)}{\sinh 2\pi v}dv.
\end{split}
\ee

\begin{lemma}[Estimate for the Laguerre coefficients]
\label{PreobLCoeffsLem}
There exists an absolute constant in the $O(1)$ such that for each nonnegative integer $k$ there exists sufficiently large $\K$ for which we have
\[
\left|\sum_{\nu=k}^{\K}{\widetilde b}_\nu(\alpha)\binom{\nu}{\nu-k}\right|
\le e^{-\alpha/2}\alpha^{O(1)}.
\]
Moreover, for $\frac{\log l}{\log T}\le\frac12$ and all sufficiently large $\K$ we can replace $s_{\K-1}\left(\alpha,(\alpha-R)\frac{\log l}{\log T}\right)$
in~\textup{\eqref{PreobLaguerreApprox}}
with its Taylor polynomial of degree $\K'\asymp\alpha^2$.
In Theorem~\ref{PreobApproximationTh}, due to the ratio of the lower incomplete gamma functions
we can truncate the sum by $\K=\alpha\left(\frac12+\epsilon\right)$.

We have the representations
\[
\begin{split}
\widetilde b_k(\alpha)=&2\sum_{\nu=0}^{k}\frac1{\nu!}\binom{k}{\nu}
\left(\Ra\right)^{\nu+1}\\
&\times\left.\frac{d^\nu}{db^\nu}\left(\frac1{b-1}e^{-\alpha/2}
{}_2F_1
\left(\left.{\scriptsize
\begin{array}{cc}
1, & 1-b \\
{} & 2-b
\end{array}
}
\right| -e^{-\alpha/2}\right)
+e^{(-\alpha/2)b}\frac{\pi}{\sin(\pi b)}
\right)\right|_{b=\Ra},
\end{split}
\]
and
\[
\begin{split}
&\sum_{\nu=k}^{\infty}{\widetilde b}_\nu(\alpha)\binom{\nu}{\nu-k}=\frac2{k!}\\
&\times\lim\limits_{r\to1{-}}\frac{d^k}{dr^k}\left(b(r)\left(\frac1{b(r)-1}e^{-\alpha/2}
{}_2F_1
\left(\left.{\scriptsize
\begin{array}{cc}
1, & 1-b(r) \\
{} & 2-b(r)
\end{array}
}
\right| -e^{-\alpha/2}\right)
+e^{(-\alpha/2)b(r)}\frac{\pi}{\sin(\pi b(r))}
\right)\right),
\end{split}
\]
where
\[
b(r)=\frac{\alpha-R}{\alpha(1-r)}.
\]

\end{lemma}

\dokvo
To prove the representations of the lemma, for $b_0=e^{(\alpha-R)/2}$ using changes of the variables
and Euler's integral representation of the ${}_2F_1$ function we write
\[
\begin{split}
\widetilde b_k(\alpha)=&\frac2{b_0}\sum_{\nu=0}^{k}\frac{(-1)^{\nu}}{\nu!}\binom{k}{\nu}
\inte_0^{b_0}(\log b_0-\log y)^\nu\frac{dy}{1+y^{\aR}}\\
=&2\sum_{\nu=0}^{k}\frac1{\nu!}\binom{k}{\nu}
\left(\Ra\right)^{\nu+1}
\left.\inte_0^1\log^\nu x x^{b-1}(1-x)^{b+1-b-1}(1-(-b_0^{\aR})x)^{-1}dx\right|_{b=\Ra}\\
=&2\sum_{\nu=0}^{k}\frac1{\nu!}\binom{k}{\nu}
\left(\Ra\right)^{\nu+1}
\left.\frac{d^\nu}{db^\nu}\left(b^{-1}
{}_2F_1
\left(\left.{\scriptsize
\begin{array}{cc}
1, & b \\
{} & b+1
\end{array}
}
\right| -e^{\alpha/2}\right)
\right)\right|_{b=\Ra}.
\end{split}
\]
Using the hypergeometric transformation formula
\[
\begin{split}
(e^{-\alpha/2})^{-b}{}_2F_1
\left(\left.{\scriptsize
\begin{array}{cc}
1, & b \\
{} & b+1
\end{array}
}
\right| -e^{\alpha/2}\right)
=&\frac{\Gamma(b-1)\Gamma(b+1)}{\Gamma^2(b)}e^{(-\alpha/2)(1-b)}
{}_2F_1
\left(\left.{\scriptsize
\begin{array}{cc}
1, & 1-b \\
{} & 2-b
\end{array}
}
\right| -e^{-\alpha/2}\right)\\
&+\frac{\Gamma(1-b)\Gamma(1+b)}{\Gamma^2(1)}e^0
{}_2F_1
\left(\left.{\scriptsize
\begin{array}{cc}
b, & 0 \\
{} & b
\end{array}
}
\right| -e^{-\alpha/2}\right)
\end{split}
\]
we obtain the stated representation of $b_k(\alpha)$.
To obtain the representation of
\[
\sum_{\nu=k}^{\infty}{\widetilde b}_\nu(\alpha)\binom{\nu}{\nu-k},
\]
we recall the generating function for the Laguerre polynomials $L_n(x)$,
\[
\frac{e^{-x/(1-r)}}{1-r}=e^{-x}\sum_{n=0}^{\infty}L_n(x)r^n.
\]
Then
\[
\begin{split}
&\sum_{\nu=k}^{\infty}{\widetilde b}_\nu(\alpha)\nu(\nu-1)\dots(\nu-k+1)r^{\nu-k}\\
&=\sum_{\nu=k}^{\infty}\left(
\inte_0^{+\infty}\left(1+\tanh\left(\frac{\alpha}{2(\alpha-R)}x-\frac{\alpha}4\right)\right)L_{\nu}(x)e^{-x}dx
\right)
\nu(\nu-1)\dots(\nu-k+1)r^{\nu-k}\\
&=2e^{-\alpha/2}\frac{d^k}{dr^k}\left(
-\inte_0^{+\infty}\frac{e^{-\frac{\alpha(1-r)(-x)}{(\alpha-R)(1-r)}}e^{-\frac x{1-r}}}
{e^{-\alpha/2}e^{-\frac{\alpha(1-r)(-x)}{(\alpha-R)(1-r)}}+1}d\left(-\frac x{1-r}\right)
\right)\\
&=2\frac{d^k}{dr^k}\left(
\frac{\alpha-R}{\alpha(1-r)}\inte_0^1
\frac{\frac{\alpha(1-r)}{(\alpha-R)}t^{1-\frac{\alpha(1-r)}{\alpha-R}}t^{\frac{\alpha(1-r)}{\alpha-R}-1}}
{1+e^{\alpha/2}t^{\frac{\alpha(1-r)}{\alpha-R}}}dt
\right),
\end{split}
\]
and the stated formula follows as above
using Euler's integral representation of ${}_2F_1$.

To obtain the estimate of the lemma,
we apply Cauchy's integral formula with a suitable contour to the established representation of the sum
and bound the Taylor coefficients by $e^{-\alpha/2}\alpha^{O(1)}$.

\section{Employing the translation}\label{PreobFEqSec}

Thus, from the functional equation~\eqref{PreobS1} we have
\[
\begin{split}
G(s)&=\sum_{l\le T}\left(\frac12+q\left(\frac{\log l}{\log T}+\delta_0(s)\right)\right)l^{-s}+O(T^{-\frac14})\\
&=e^{\alpha/2}\sum_{l\le T}\left(\frac12-\widetilde q\left(\frac{\log l}{\log T}+\delta_1(s)\right)\right)l^{-(s+\ds)}\\
&+\RR,
\end{split}
\]
where
\[
\RR=\sum_{l\le T}{\mathcal R}_{K,\K,\alpha}\left(\frac{\log l}{\log T}-\frac12\right)l^{-s}
\]
and
\[
{\mathcal R}_{K,\K,\alpha}\ll e^{-\alpha}\quad\text{for some}\quad K\asymp\alpha,\K\asymp\alpha^3.
\]

We denote
\[
G^{*}(s+\ds)=e^{\alpha/2}\sum_{l\le T}\left(\frac12-\widetilde q\left(\frac{\log l}{\log T}\right)\right)l^{-(s+\ds)}.
\]

Theorem~\ref{PreobTranslationTh}, Subsection~\ref{PreobOutFEq} asserts that the terms $\delta_1(s)$ and $\RR$ do not affect
the principal inequality of the Levinson--Conrey method, as per~\cite[Chapter 18, (18{.}14)--(18{.}19)]{PreobIw14}.

For $\ds=\frac{\alpha}{\log T}$ we can write
\[
\frac{G^{*}(s+\ds)+O(T^{-\frac14})}{\zeta(s+\ds)}=c_0(\alpha)+\lambda(s+\ds),
\]
where
\[
\lambda(s+\ds)
=\frac{c_1(\alpha)}{\log T}\frac{\zeta'}{\zeta}(s+\ds)
+\frac{c_2(\alpha)}{(\log T)^2}\frac{\zeta''}{\zeta}(s+\ds)+\dots+\frac{c_{\K}(\alpha)}{(\log T)^{\K}}\frac{\zeta^{(\K)}}{\zeta}(s+\ds)
\]
and
\[
c_k(\alpha)=\frac{e^{\alpha/2}}{2}(\alpha-R)^k\frac1{k!}\sum_{\nu=k}^{\K}{\widetilde b}_\nu(\alpha)\binom{\nu}{\nu-k}
\]
with $\K\asymp\alpha^3$.


We now employ a generalization of Selberg's construction~\cite{PreobSel46}
to approximate the function $\lambda(s+\ds)$ by a Dirichlet polynomial on a set which has the measure at least $(1-\alpha^{1-C(1-\epsilon)})T$,
where $\alpha>0$ grows with $T$ sufficiently slowly.

\section{The Selberg approximation}\label{selbergsec}

We have
\[
\frac{\zeta''}{\zeta}=\left(\frac{\zeta'}{\zeta}\right)^{(1)}+\left(\frac{\zeta'}{\zeta}\right)^2,
\]
\[
\frac{\zeta^{(3)}}{\zeta}=\left(\frac{\zeta'}{\zeta}\right)^{(2)}+3\left(\frac{\zeta'}{\zeta}\right)^{(1)}\left(\frac{\zeta'}{\zeta}\right)
+\left(\frac{\zeta'}{\zeta}\right)^3.
\]
In general, by Fa\`a di Bruno's formula we have
\[
\frac1{\zeta(s)}\left(\exp(\log\zeta(s))\right)^{(k)}=
\sum_{\substack{R_1\ge0,\ldots,R_k\ge0\\R_1+\dots+kR_k=k}}
\frac{k!}{R_1!\dots R_k!}\prod_{j=1}^k\left(\frac{\left(\frac{\zeta'}{\zeta}\right)^{(j-1)}(s)}
{j!}\right)^{R_j}.
\]
Let
\[
\left(\frac{\zeta'}{\zeta}\right)^{\wedge(R_1,\ldots,R_k)}
\]
denote the product of the zeroth derivative of $\frac{\zeta'}{\zeta}$ to the power $R_1$,
the $1$st derivative to the power $R_2$, $\ldots$, the $(k-1)$th derivative to the power $R_k$.
Applying Fa\`a di Bruno's formula, we obtain the expression with the coefficients
$c^{(1)}_{R_1,\ldots,R_k}=\frac{k!}{R_1!{1!}^{R_1}\dots R_k!{k!}^{R_k}}$:
\[
\frac{\zeta^{(k)}}{\zeta}=\sum_{\substack{R_1\ge0,\ldots,R_k\ge0\\R_1+\dots+kR_k=k}}
c^{(1)}_{R_1,\ldots,R_k}\left(\frac{\zeta'}{\zeta}\right)^{\wedge(R_1,\ldots,R_k)},
\]
which yields the expression for $\lambda(s+\ds)$ of the form
\be
\label{Preoblex}
\begin{split}
\lambda(s+\ds)
&=\frac{c_1(\alpha)}{\log T}\frac{\zeta'}{\zeta}(s+\ds)
+\frac{c_2(\alpha)}{(\log T)^2}\frac{\zeta''}{\zeta}(s+\ds)+\dots+\frac{c_{\K}(\alpha)}{(\log T)^{\K}}\frac{\zeta^{(\K)}}{\zeta}(s+\ds)\\
&=\sum_{k=1}^{\K}\sum_{\substack{R_1\ge0,\ldots,R_k\ge0\\R_1+\dots+kR_k=k}}
\frac{c^{(2)}_{k;R_1,\ldots,R_k}}{(\log T)^k}\left(\frac{\zeta'}{\zeta}\right)^{\wedge(R_1,\ldots,R_k)}
\end{split}
\ee
with the coefficients $c^{(2)}_{k;R_1,\ldots,R_k}$.

We now give the well-known Selberg formula
for a Dirichlet polynomial approximation of the function $\zeta'/\zeta(\sigma+it)$
and a lemma on the measure of the set of $t$ for which the approximation can fail.
We then generalize these results to $\lambda(s+\ds)$.

First we define
\[
\sigma_{x,t}=\frac12+2\max\limits_{\rho^{\ast}}\left(\beta^{\ast}-\frac12,\frac2{\log x}\right),
\]
where $x\ge2$ and $t>0$. The maximum is taken over all such zeros $\rho^{\ast}$ of the zeta-function that satisfy $|\gamma^{\ast}-t|\le x^{3|\beta^{\ast}-1/2|}/\log x$.
\begin{lemma}[A. Selberg~\cite{PreobSel46}]\label{PreobSelSimple}
If $\sigma_{x,t}\le\sigma$ and $2\le x\le t^2$, then
\[
\begin{split}
-\frac{\zeta'}{\zeta}(\sigma+it)=\sum_{n\le x^3}\frac{\Lambda_x(n)}{n^{\sigma+it}}
&+O\left(x^{(1/2-\sigma)/2}\left|\sum_{n\le x^3}\frac{\Lambda_x(n)}{n^{\sigma_{x,t}+it}}\right|\right)\\
&+O\left(x^{(1/2-\sigma)/2}\log t\right),
\end{split}
\]
where
\[
\Lambda_x(n)=\begin{cases}
\Lambda(n)&\text{if $n\le x$},\\
\Lambda(n)\frac{\log^2(x^3/n)-2\log^2(x^2/n)}{2\log^2 x}&\text{if $x<n\le x^2$},\\
\Lambda(n)\frac{\log^2(x^3/n)}{2\log^2 x}&\text{if $x^2<n\le x^3$}.
\end{cases}
\]
\end{lemma}

\begin{lemma}[Selberg--Jutila zero-density estimate~\cite{PreobJut82}]
\label{PreobSelJutLem}
\[
N(\sigma,T)\ll T^{1-(1-\epsilon)(\sigma-1/2)}\log T.
\]
\end{lemma}

We now give the following lemma on the measure of the set of $t\in(2,T)$ for which $\sigma_{x,t}>\sigma$
(see~\cite[Lemma 2{.}4]{PreobLes13b}).
The proof of the lemma uses the Selberg--Jutila zero-density estimate~\cite{PreobJut82}.
\begin{lemma}[S. Lester]
\label{PreobLesLemma}
Let $1/2+4/(\log x)\le\sigma\le2$, and for a fixed $0<\epsilon<1$
let $10\le x\le T^{\epsilon/3}$.
Then
\[
\meas\{t\in(2,T)\ {:}\ \sigma_{x,t}>\sigma\}\ll T^{1-(1/2-\epsilon)(\sigma-1/2)}\frac{\log T}{\log x}
\]
with the implied constant depending only on $\epsilon$.
\end{lemma}

\begin{lemma}[Tsang~\cite{PreobTsang84}, pp.~68--69 (Lemma 5.4)]
For $t\in[T,2T]$, $x=T^{1/\alpha^{1+c}}$,
$t\ne\gamma$ we have
\[
\log\zeta(\sigma_1+it)=\begin{cases}
&\sum_{n\le x^3}\frac{\Lambda_x(n)}{n^{\sigma_1+it}\log n}
+O\left(\frac{x^{\frac12\left(\frac12-\sigma_1-it\right)}}{\log x}
\left|\sum_{n\le x^3}\frac{\Lambda_x(n)}{n^{\sigma_1+it}}\right|\right)\\
&+O\left(\frac{x^{\frac12\left(\frac12-\sigma_1-it\right)}}{\log x}\log(\sigma_1+it)\right)
\quad\text{for}\quad\sigma_1>\sigma_{x,t},\\
&\sum_{n\le x^3}\frac{\Lambda_x(n)}{n^{\sigma_{x,t}+it}\log n}
+O\left(\left(\sigma_{x,t}-\sigma_1\right)\left|\sum_{n\le x^3}\frac{\Lambda_x(n)}{n^{\sigma_{x,t}+it}}\right|\right)\\
&+O\left(\inte_{\sigma_1+it}^{\sigma_{x,t}+it}\log s\,ds\right)\\
&-\sum_{\rho}\inte_{\sigma_1}^{\sigma_{x,t}}\frac{\sigma_{x,t}-u}{(u+it-\rho)(\sigma_{x,t}+it-\rho)}du
\quad\text{for}\quad\frac12\le\sigma_1<\sigma_{x,t}.
\end{cases}
\]
\end{lemma}

Next, using Tsang's lemma, the Selberg--Jutila zero-density estimate and the higher order Cauchy integral formula
we deduce Theorem~\ref{PreobApproximationTh}.

\textbf{Proof of Theorem~\ref{PreobApproximationTh}.}
First we write
\[
\log\zeta\left(\frac12+\frac{\alpha-R}{\log T}+it\right)
=\log\zeta\left(\frac12+\frac{\alpha^{1/2+4C}}{\log T}+it\right)
+\log\zeta\left(\frac12+\frac{\alpha-R}{\log T}+it\right)-\log\zeta\left(\frac12+\frac{\alpha^{1/2+4C}}{\log T}+it\right).
\]
The derivatives of the first term are well approximated by the Dirichlet polynomial of length $x_0=T^{1/\alpha^{4C}}$.
Then we apply the higher order Cauchy integral formula
to the difference $\log\zeta(\sigma_1+it_1)-\log\zeta\left(\sigma_1+\frac{\alpha^{1/2+4C}-\alpha+R}{\log T}+it_1\right)$ in the disc
\[
\left\{s_1=\sigma_1+it_1\ :\ %
\left|s_1-\left(\frac12+\frac{\alpha-R}{\log T}+it\right)\right|\le\frac{\alpha-R-C\log\alpha}{\log T}\right\},
\]
if there are no zeros $\rho$ in the disc. The Selberg--Jutila zero-density estimate
implies that this can fail for at most $\frac1{\alpha^{C(1-\epsilon)}}T\log T$ disjoint discs
and hence for the set of $t\in[T,2T]$ which has the measure at most $\frac1{\alpha^{C(1-\epsilon)-1}}T$.
Cauchy's formula and Tsang's lemma with $x=T^{1/\alpha^{1/2+4C}}$ give
\[
\begin{split}
&\left(\log\zeta(\sigma+\ds+it)-\log\zeta\left(\sigma+\ds+\frac{\alpha^{1/2+4C}-\alpha+R}{\log T}+it\right)\right)^{(j)}\\
&\ll\left(\sigma_{x,t_1}-\frac12\right)\frac{j!(\log T)^j}{(\alpha-R-C\log\alpha)^j}
\left|\sum_{n\le x^3}\frac{\Lambda_x(n)}{n^{\sigma_{x,t_1}+it}}\right|\\
&+\left|\left.\left(\frac{d}{ds_1}\right)^{j}
\inte_{s_1}^{s_1+\frac{\alpha^{1/2+4C}-\alpha+R}{\log T}}\log s\,ds\right|_{s_1=\sigma_1+it_1}\right|\\
&+\frac{j!(\log T)^j}{(\alpha-R-C\log\alpha)^j}
\left|\sum_{\rho}\inte_{\sigma_1}^{\sigma_1+\frac{\alpha^{1/2+4C}-\alpha+R}{\log T}}\frac{\sigma_{x,t_1}-u}{(u+it_1-\rho)(\sigma_{x,t_1}+it_1-\rho)}du\right|
\end{split}
\]
for some real numbers
$\sigma_1\in\left[\frac12+\frac{C\log\alpha}{\log T},\frac12+\frac{2(\alpha-R)-C\log\alpha}{\log T}\right]$
and $t_1\in\left[t-\frac{\alpha-R-C\log\alpha}{\log T},t+\frac{\alpha-R-C\log\alpha}{\log T}\right]$.
The term with $\rho$ is difficult to analyze directly. However, if we consider the approximation to
\be
\label{PreobSecDiff}
\log\zeta\left(\frac12+\frac{\alpha^{1/2+4C}}{\log T}+it_1\right)-\log\zeta\left(\frac12+\frac{\alpha^{1/2+8C}}{\log T}+it_1\right)
\ee
using the value $x_1=T^{1/\alpha^{1/2+8C}}$
then the moments of the latter terms with $\rho$ and $\log s$ are seen to be of the same order of magnitude as the moments of the former terms with $\rho$ and $\log s$.
Applying Tsang's lemma again to the difference~\eqref{PreobSecDiff}, but with the value $x_0=T^{1/\alpha^{4C}}$,
we get
\[
\begin{split}
&\left|\sum_{\rho}\inte_{\sigma_1}^{\sigma_1+\frac{\alpha^{1/2+4C}-\alpha+R}{\log T}}\frac{\sigma_{x,t_1}-u}{(u+it_1-\rho)(\sigma_{x,t_1}+it_1-\rho)}du\right|\\
&\ll\left|\sum_{n\le x_0^3}\frac{\Lambda_{x_0}(n)}{\log n}
\left(\frac1{n^{1/2+\alpha^{1/2+4C}/\log T+it_1}}-\frac1{n^{1/2+\alpha^{1/2+8C}/\log T+it_1}}\right)\right|\\
&+\frac{\alpha^{1/2+4C}}{\log T}\left|\sum_{n\le x_0^3}\frac{\Lambda_{x_0}(n)}{n^{1/2+\alpha^{1/2+4C}/\log T+it_1}}\right|
+\frac{\alpha^{1/2+4C}}{\log T}\left|\sum_{n\le x_0^3}\frac{\Lambda_{x_0}(n)}{n^{1/2+\alpha^{1/2+8C}/\log T+it_1}}\right|\\
&+\frac{\alpha^{1/2+8C}}{\log T}\left|\sum_{n\le x_1^3}\frac{\Lambda_{x_1}(n)}{n^{1/2+\alpha^{1/2+8C}/\log T+it_1}}\right|\\
&+e^{-\sqrt{\alpha}}\left|
\inte_{\frac12+\frac{\alpha^{1/2+4C}}{\log T}+it_1}^{\frac12+\frac{\alpha^{1/2+8C}}{\log T}+it_1}\log s\,ds\right|
\end{split}
\]
for a set of $t\in[T,2T]$ which has the measure at least $T(1-\alpha^{1-C(1-\epsilon)})$.

Here we applied Lemma~\ref{PreobLesLemma} to control the values $\sigma_{x,t_1}$.
Now Theorem~\ref{PreobApproximationTh} follows using Taylor expansion of the Dirichlet polynomials in powers of $(t-t_1)$ and the fact that
\emph{we can multiply the short Dirichlet polynomials by suitable quantities} so that the absolute value of the sum of the multiplied polynomials
is larger than the sum of the absolute values of the original polynomials for almost all $t$.\qed

Next we formulate the following property of exceptional $t$'s
for which the values of the approximating Dirichlet polynomial
\[
A(\sigma+\ds+it)
\]
defined in Section~\ref{PreobOutSeLester} lie exterior to an appropriate closed Jordan region ${\mathcal C}\cup{\mathcal C}'\cup{\mathcal C}''\subset{\mathbb C}$
not containing the point $-c_0(\alpha)$,
and/or for which the approximation in Theorem~\ref{PreobApproximationTh} fails.
Theorem~\ref{PreobExceptionTh} is motivated by results in~\cite[Chapter 2]{PreobLes13b},
which provide Gaussian distribution for
\[
\frac{\zeta'}{\zeta}(\sigma+it)
\]
and the Dirichlet polynomial
\[
\sum_{n\le X}\frac{\Lambda_x(n)}{n^{\sigma+it}}.
\]
Lester's results are valid for larger values of $\sigma$, though.

\begin{figure}
 \begin{center}
   \includegraphics[width=0.6\textwidth,keepaspectratio]{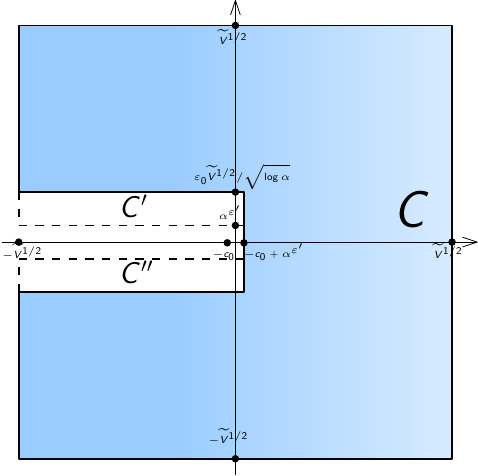}
   \caption{The sets ${\mathcal C}'$, ${\mathcal C}''$ and ${\mathcal C}$} \label{PreobC0CFig}
 \end{center}
\end{figure}

\begin{theo}
\label{PreobExceptionTh}
For a fixed small $\epsilon>0$ let $R=\epsilon\log\alpha$
and let $\alpha$ go to infinity with $T$ sufficiently slowly.
Let the closed Jordan region ${\mathcal C}\subset{\mathbb C}$
not containing the point $-c_0(\alpha)=-\frac{e^{\alpha/2}}{2}\left(1-\tanh\left(\frac{\alpha}4\right)\right)$
be the square
\[
U=\{z\in{\mathbb C}\ :\ -\widetilde{V}^{1/2}\le\Re z\le\widetilde{V}^{1/2}\text{ and }-\widetilde{V}^{1/2}\le\Im z\le \widetilde{V}^{1/2}\}
\]
from which we remove the box $S_{\eps_0}$ of the form
\[
\{z\in{\mathbb C}\ :\ -\widetilde{V}^{1/2}\le\Re z<-c_0(\alpha)+\alpha^{\epsilon'}\text{ and }-\eps_0\widetilde{V}^{1/2}/\sqrt{\log\alpha}<\Im z<\eps_0\widetilde{V}^{1/2}/\sqrt{\log\alpha}\}
\]
where $\widetilde{V}=\widetilde{V}(\alpha)>0$, $\eps_0=\eps_0(\alpha)>0$, that is,
\[
{\mathcal C}=U\setminus S_{\eps_0}.
\]
See Figure~\ref{PreobC0CFig}.

Let
\[
  {\mathcal C}'=\{z\in{\mathbb C}\ :\ -\widetilde{V}^{1/2}\le\Re z\le -c_0(\alpha)+\alpha^{\epsilon'}\text{ and }\alpha^{\eps'}\le\Im z\le\eps_0\widetilde{V}^{1/2}/\sqrt{\log\alpha}\}
\]
and
\[
  {\mathcal C}''=\{z\in{\mathbb C}\ :\ -\widetilde{V}^{1/2}\le\Re z\le -c_0(\alpha)+\alpha^{\epsilon'}\text{ and }-\eps_0\widetilde{V}^{1/2}/\sqrt{\log\alpha}\le\Im z\le -\alpha^{\eps'}\}
\]

For $\delta_T\ll1$ we define the exceptional set $E_{\alpha,\delta_T}$
to be such that for $\sigma+\ds=\frac12+\frac{\alpha-R+\delta_T}{\log T}$ we have
\[
A(\sigma+\ds+it)\in{\mathcal C}\subset{\mathbb C}
\]
if and only if $t$ is in the set $[T,2T]\setminus E_{\alpha,\delta_T}$.

Let the sequence of polynomials ${\mathcal L}_{n,{\mathcal C}}(w)$
defined in Section~\textup{\ref{PreobRungeSection}} uniformly minorizes
\[
f(w)=\frac{\widetilde{\mathcal M}}{c_0(\alpha)+w}
\]
in ${\mathcal C}$, where $\widetilde{\mathcal M}=\alpha^{\epsilon'}$
and
\[
  \max_{w\in{\mathcal C}'\cup{\mathcal C}''}|{\mathcal L}_{n,{\mathcal C}}(w)|\le 1.
\]
Then there exist $\widetilde{V}=\alpha^{C(1/2+2\epsilon)+c}$ with some fixed $c=o(C)>0$ and $\eps_0=\eps_0(\alpha)\to0$, $\eps'=\eps'(\alpha)\to0$ as $\alpha\to +\infty$
such that for our choice of the function $F^{\ast}$
the terms $\frac2{R}\eps_E\log I_E(R)$ and $\frac2{R}L_E(R)$ in Theorem~\textup{\ref{PreobLevCon}}
for the exceptional set $E=E_{\alpha,\delta_T}\cup([T,2T]\setminus{\mathcal M}_\alpha)$,
with ${\mathcal M}_\alpha$ as in Theorem~\textup{\ref{PreobApproximationTh}}, are at most $\epsilon(\alpha,T)$,
where $\epsilon(\alpha,T)$ can be made arbitrarily small.
\end{theo}

In the proof of Theorem~\ref{PreobExceptionTh},
we shall apply a method of proof of the following results of Lester~\cite[Theorems 1 and 2]{PreobLes13a}:

\begin{theo}
\label{PreobLesterThRect}
Let $\psi(T)=(2\sigma-1)\log T$, and for $\psi(T)\geq 1$ define
\[
V=V(\sigma)=\displaystyle\frac12\sum_{n=2}^{\infty}\frac{\Lambda^2(n)}{n^{2\sigma}},
\]
\[
\bOmega=e^{-10}\min\big(V^{3/2},(\psi(T)/\log\psi(T))^{1/2}\big).
\]
Suppose that $\psi(T)\to\infty$ with $T$,  $\psi(T)=o(\log T)$,
and that $R$ is a rectangle in $\mathbb C$ whose sides are parallel to the coordinate axes.
Then we have
\be
\label{PreobRectMeas}
\begin{split}
\meas\bigg\{t\in(0,T) : \frac{\zeta'}{\zeta}(\sigma+it)V^{-1/2}\in R\bigg\}
=\frac{T}{2\pi}\iint_R\! e^{-(x^2+y^2)/2}\, dx\, dy+O\left(T\frac{(\meas(R)+1)}{\bOmega}\right).
\end{split}
\ee
\end{theo}

\begin{theo}
\label{PreobLesterThDisk}
Let $\psi(T)=(2\sigma-1)\log T$, and for $\psi(T)\geq 1$ define
\[
V=V(\sigma)=\displaystyle\frac12\sum_{n=2}^{\infty}\frac{\Lambda^2(n)}{n^{2\sigma}},
\]
\[
\bOmega=e^{-10}\min\big(V^{3/2},(\psi(T)/\log\psi(T))^{1/2}\big).
\]
Suppose that $\psi(T)\to\infty$ with $T$,  $\psi(T)=o(\log T)$,
and that $r$ is a real number such that $r\bOmega\geq 1$.
Then we have
\be
\label{PreobMainMeas}
\begin{split}
\meas\bigg\{t\in(0,T) : \bigg|\frac{\zeta'}{\zeta}(\sigma+it)\bigg|\leq\sqrt Vr\bigg\}
=\, T(1-e^{-r^2/2})+O\left(T\left(\frac{r^2+r}{\bOmega}\right)\right).
\end{split}
\ee
If, in addition, we let $\widetilde\bOmega=
\min\big((2\sigma-1)e^{\sigma/(2\sigma-1)},e^{-10}(\psi(T)/\log\psi(T))^{1/2}\big)$,
then we have for $r\widetilde\bOmega\geq 1$
\be
\label{PreobSmallMeas}
\meas\bigg\{t\in[0,T] : \bigg|\frac{\zeta'}{\zeta}(\sigma+it)\bigg|\leq\sqrt Vr\bigg\}
\ll Tr^2.
\ee
\end{theo}

\textbf{Proof of Theorem~\textup{\ref{PreobExceptionTh}}.}
We consider the distribution of values of the Dirichlet polynomial
\[
A(\sigma+\ds+it)
\]
defined in Subsection~\ref{PreobOutSeLester}.
We emphasize that in the construction of $A(\sigma+\ds+it)$ we use {\it very short} Dirichlet polynomials
for $\left(\frac{\zeta'}{\zeta}\right)^{(j)}(s+\Delta\sigma)$ in multivariate polynomial~\eqref{Preoblex}
for $\lambda(s+\Delta\sigma)$.
The lengths of those short Dirichlet polynomials are as small as $T^{1/N}$
with $N\asymp\alpha^{1+c}$, where $c>0$ is fixed.
It is important to remark that with this choice of the lengths,
$A(\sigma+\ds+it)$ will {\it not} be a precise approximation to $\lambda(\sigma+\Delta\sigma+it)$.
However, Theorem~\ref{PreobApproximationTh}
implies that for almost all $t\in[T,2T]$ we have
\be
\label{PreobSelbergUpperBound}
|\lambda(\sigma+\Delta\sigma+it)|\ll|A(\sigma+\ds+it)|
\ee
with the implied constant being absolute.
This is a crucial fact about $A(\sigma+\ds+it)$ that will be used in Section~\ref{seccompletion}.

For any large but fixed $C>0$ from the statement of Theorem~\ref{PreobApproximationTh} there exists $c=o(C)$
such that for $\widetilde{V}=\alpha^{C(1/2+2\epsilon)+c}$ we have the rough estimates, analogous to~\eqref{PreobMainMeas} and~\eqref{PreobRectMeas}, respectively:
\be
\label{PreobAExtDisk}
\meas\{t\in[T,2T]\ :\ |A(\sigma+\ds+it)|>\widetilde{V}^{1/2}\}
\ll T\alpha^{-C(1-\epsilon)},
\ee
and
\be
\label{PreobANarrowStrip}
\meas\{t\in[T,2T]\ :\ A(\sigma+\ds+it)\in S_{\eps_0}\setminus({\mathcal C}'\cup{\mathcal C}'')\}
\ll\frac{\alpha^{\eps'}}{\widetilde{V}^{1/2}} T
\ee
(see Figure~\ref{PreobC0CFig}).
To prove the existence of $c=o(C)>0$, $\widetilde{V}=\alpha^{C(1/2+2\epsilon)+c}$, $\eps_0=\eps_0(\alpha)\to0$ and $\epsilon'(\alpha)\to0$ as $\alpha\to +\infty$
we use the formulas in Sections~\ref{PreobFEqSec}, \ref{selbergsec},
Theorem~\ref{PreobApproximationTh}
and a generalization of~\cite[Section 2{.}3{.}1, Lemma 2{.}5]{PreobLes13b} to study the distribution of the short Dirichlet polynomials.

In our function $G^{{*}}(s)$
defined in the beginning of Section~\ref{PreobFEqSec}
we can take $\K=\alpha(1/2+\epsilon)$ in the polynomial $\widetilde{q}$.

Thus we have
\[
\begin{split}
&\meas\bigg\{t\in[T,2T] : \left|\frac{A_j(\sigma+\ds+it)}{j!c\sqrt{\log\alpha}}\right|^{R_j}
>V^{R_j/2}\bigg\}\\
&\ll T\alpha^{-C(1-\epsilon)}.
\end{split}
\]
Bounding the Laguerre coefficients using Lemma~\ref{PreobLCoeffsLem},
and using the identity
\[
\sum_{\substack{R_1\ge0,\ldots,R_k\ge0\\R_1+\dots+kR_k=k}}
\frac{(R_1+\dots+R_k)!}{R_1!\dots R_k!}=2^{k-1},
\]
we deduce that the terms of the formula in Section~\ref{PreobFEqSec} contribute $\frac{\widetilde{V}^{1/2}}{\alpha^{O(1)}}$ with the exceptional measure $\alpha^{-C(1-\epsilon)}$.

Then from~\eqref{PreobAExtDisk} and~\eqref{PreobANarrowStrip}
we get the following bound on the measure of the exceptional set:
\[
\begin{split}
&\meas E_{\alpha,\delta_T}=\meas\{t\in(T,2T)\ \colon\ A(\sigma+\ds+it)\not\in{\mathcal C}\}\\
&\ll\eps_0T+\frac{\alpha^{\eps'}}{\widetilde{V}^{1/2}}T+T\alpha^{-C(1-\epsilon)}.
\end{split}
\]
For the purpose of applying the Cauchy-Schwarz inequality, we compute the mean square $I$
of the mollifier $M(s)$ which is the optimal one
for $\zeta(s)$ on the line $\Re s=\frac12+\frac{\alpha-R}{\log T}$,
times $G_{\alpha}(s)$ on the line $\Re s=\frac12-\frac{R}{\log T}$,
where
\[
G_{\alpha}(s)=\zeta(s)-\sum_{l=1}^{\infty}\tanh\left(\frac\alpha2\left(\frac{\log l}{\log T}-\frac12\right)\right)l^{-s}.
\]
Here we use the integrand without the shift, though this integrand is approximately equal to the shifted one
by the translation functional equation and Theorem~\ref{PreobTranslationTh}.

The optimal function $P(x)$ for $\zeta(s)$ on the line $\Re s=\frac12+\frac{\alpha-R}{\log T}$
is given by (see~\cite[Theorem 2]{PreobCon89})
\[
P(x)=\frac{\sinh(\alpha-R)x\theta}{\sinh(\alpha-R)\theta}.
\]
Let
\[
Q(y)=\frac{1-\tanh\left(\frac\alpha2\left(y-\frac12\right)\right)}{1+\tanh\frac\alpha4}
\]
and
\beq
\label{PreobMeanSqInt}
A=\int_0^1(Q(y)e^{Ry})^2dy,\qquad A_1=\int_0^1\bigg((Q(y)e^{Ry})'\bigg)^2dy.
\eeq
Then for $\alpha$ independent of $T$ we get
\[
A\asymp e^{R}/R,\qquad A_1\asymp\alpha e^{R}+Re^R,
\]
\[
P'(x)=(\alpha-R)\theta\frac{\cosh(\alpha-R)x\theta}{\sinh(\alpha-R)\theta}
\]
and
\[
\int_0^1P'(x)^2dx=(\alpha-R)^2\theta^2\frac{(\alpha-R)\theta/(\sinh(\alpha-R)\theta)+\cosh(\alpha-R)\theta}{2(\alpha-R)\theta\sinh(\alpha-R)\theta}\asymp\alpha.
\]
Thus
\[
I=1+A\int_0^1P'(x)^2dx+A_1\int_0^1P(x)^2dx+(Q(1)e^R-Q(0))(P(1)-P(0))\asymp\alpha e^{R}.
\]
Hence the term $\frac2{R}\eps_E\log I_E(R)$ in Theorem~\textup{\ref{PreobLevCon}}
is such that
\be
\label{PreobI}
\frac2{R}\eps_E\log I_E(R)\ll R^{-1}
\left(\eps_0+\frac{\alpha^{\eps'}}{\widetilde{V}^{1/2}}+\alpha^{-C(1-\epsilon)}\right)\log\left(\alpha e^{R}\right).
\ee

To bound the term $\frac2{R}L_E(R)$ in Theorem~\textup{\ref{PreobLevCon}}
we need to know the degree $n$
of the approximating polynomial ${\mathcal L}_{n,{\mathcal C}}(w)$
giving an acceptable error,
and a bound on the coefficients of this polynomial.
By Cauchy's integral formula,
for ${\mathcal M}$ being the maximum of the absolute values of the coefficients we have
\[
{\mathcal M}=O\left(\max_{w\in S_{\eps_0}\setminus({\mathcal C}'\cup{\mathcal C}'')}|{\mathcal L}_{n,{\mathcal C}}(w)|\right),
\]
since for $w\in{\mathcal C}\cup{\mathcal C}'\cup{\mathcal C}''$ (see Figure~\ref{PreobCFig}) ${\mathcal L}_{n,{\mathcal C}}(w)$ is close to the function $f(w)=\frac{\widetilde{\mathcal M}}{c_0(\alpha)+w}$.

Using our estimate~\eqref{PreobUniNorm} in Section~\ref{PreobRungeSection} below,
we choose the degree $n$ of the polynomial ${\mathcal L}_{n,{\mathcal C}}(w)$:
\[
n\asymp\frac{\widetilde{V}^{1/2}}{\alpha^{\eps'}}e^{\sqrt{\log\alpha}/\eps_0}\log\left(\widetilde{V}\right).
\]
Next, we define
\[
E_S=\{t\in(T,2T)\ \colon\ A(\sigma+\ds+it)\in S_{\eps_0}\setminus({\mathcal C}'\cup{\mathcal C}'')
\text{ and }\log|{\mathcal L}(A(\sigma+\ds+it))|>0\}
\]
and estimate $L_E(R)$ as follows
\[
\begin{split}
L_E(R)\ll\frac1{T}&\left(\meas E_S\max_{z\in S_{\eps_0}\setminus({\mathcal C}'\cup{\mathcal C}'')}\log|{\mathcal L}_{n,{\mathcal C}}(z)|+(\meas E\setminus E_S)\log{\mathcal M}\right.\\
&\left.+(\meas E\setminus E_S)^{1-1/k}n\left(\inte_T^{2T}\log^k|A(\sigma+\ds+it)|\,dt\right)^{1/k}\right)
\end{split}
\]
for any $k\ge2$.
As above, from~\eqref{PreobAExtDisk} and~\eqref{PreobANarrowStrip}, and using Jensen's inequality we get
\be
\label{PreobL}
\begin{split}
L_E(R)\ll\frac1{T}&\bigg(\left(\frac{\alpha^{\eps'}}{\widetilde{V}^{1/2}}+\alpha^{-C(1-\epsilon)}\right)T\max_{z\in S_{\eps_0}\setminus({\mathcal C}'\cup{\mathcal C}'')}\log|{\mathcal L}_{n,{\mathcal C}}(z)|\\
&+T^{1-1/k}(\alpha^{-C(1-\epsilon)})^{1-1/k}nT^{1/k}\log\left(\widetilde{V}\right)\bigg).
\end{split}
\ee
The bound
\[
\max_{z\in S_{\eps_0}\setminus({\mathcal C}'\cup{\mathcal C}'')}\log|{\mathcal L}_{n,{\mathcal C}}(z)|\ll\log\widetilde{\mathcal M}+n\log\left(1+\frac{\alpha^{\eps'}}{\widetilde{V}^{1/2}}\right)
\]
can be seen by re-expanding the polynomial ${\mathcal L}_{n,{\mathcal C}}(z)$
in powers of $z-z_0$ where
\[
-\frac12\widetilde{V}^{1/2}\le\Re z_0\le-c_0(\alpha)+\alpha^{\eps'}
\]
and $\Im z_0=-\frac12\widetilde{V}^{1/2}$, and bounding the derivatives at $z_0$ by Cauchy's formula
on the circle $|z-z_0|=\frac12\widetilde{V}^{1/2}-2\alpha^{\eps'}$ using the fact that
${\mathcal L}_{n,{\mathcal C}}(z)$ is close to the function $f(z)=\frac{\widetilde{\mathcal M}}{c_0(\alpha)+z}$
in~${\mathcal C}$.

Recalling
\[
\widetilde{V}=\alpha^{C(1/2+2\epsilon)+c},\quad c=o(C),
\]
Theorem~\ref{PreobExceptionTh} follows upon choosing
large $\alpha$, fixed $C>0$, $R=\epsilon\log\alpha$ with a small fixed $\epsilon>0$, $\K=\alpha(1/2+2\epsilon)$
in~\eqref{Preoblex} and~\eqref{PreobLaguerreApprox},
$\eps_0=\eps_0(\alpha)\to0$ as $\alpha\to +\infty$ in~\eqref{PreobL} and~\eqref{PreobI}, and $\epsilon'=\epsilon'(\alpha)\to0$ as $\alpha\to +\infty$.

\section{Runge's approximation polynomials}\label{PreobRungeSection}

We substitute the Dirichlet polynomial $A(\sigma+\ds+it)$ for $z$ in the sequence of polynomials
${\mathcal L}_{n,{\mathcal C}}(z)$ obtained in the following two-stage construction:
in the first stage, we construct Runge's approximation polynomials ${\mathcal R}_{n,{\mathcal C}}(w)$ for the function
\[
f(w)=\frac{\widetilde{\mathcal M}}{c_0(\alpha)+w},
\]
uniformly approximating this function in the closed Jordan region ${\mathcal C}\subset{\mathbb C}$
not containing the point $-c_0(\alpha)$.
In the second stage, we choose the translation parameter $z_0=z_0(n)=O(\widetilde{V}^{1/2})$ and scale factor $a=a(n)=O(1)$
so that
\[
  {\mathcal L}_{n,{\mathcal C}}(z)=\frac1{\max_{z: az+z_0\in {\mathcal C}'\cup{\mathcal C}''}|{\mathcal R}_{n,{\mathcal C}}(az+z_0)|}
  {\mathcal R}_{n,{\mathcal C}}(az+z_0)
\]
satisfies ${\mathcal L}_{n,{\mathcal C}}(0)=\frac1{c_0(\alpha)}$.
The parameters exist by the Hadamard three-circle theorem.

The classical problem of Runge's approximation can be solved explicitly using Lagrange's interpolation of $f$
and loci $C_R$ of Green's function for our polygonal region ${\mathcal C}$.

The error is estimated in terms of the increment $\Delta R=R-1$ where $R>1$ is a value
for which the singularity $-c_0(\alpha)$ of $f$ does not lie on or within $C_R$.

Using the following theorems~\cite[\S\,4{.}1, Theorem~1 and \S\,4{.}5, Theorems~4, 5]{PreobWal56}
and~\cite[Ch.\,II, \S\,3, Theorem~1, Steps~1, 2]{PreobGai87}
we supply a bound on the error in this approximation.

\begin{theo}
\label{PreobWalshTh1}
Let ${\mathcal C}$ be a closed limited point set of the $z(=x+iy)$-plane whose
complement $K$ \textup{(}with respect to the extended plane\textup{)} is connected and regular in the
sense that $K$ possesses a Green's function $G(x,y)$ with pole at infinity. Then the
function $w=\phi(z)=e^{G+iH}$, where $H$ is conjugate to $G$ in $K$, maps $K$ conformally but
not necessarily uniformly onto the exterior of the unit circle $\gamma$ in the $w$-plane so that
the points at infinity in the two planes correspond to each other\textup{;} interior points of
$K$ correspond to exterior points of $\gamma$, and exterior points of $\gamma$ correspond to interior
points of $K$.

Each equipotential locus in $K$ such as $C_R\colon G(x,y)=\log R>0$, or $|\phi(z)|=R>1$,
either consists of a finite number of finite mutually exterior analytic Jordan
curves or consists of a finite number of contours which are mutually exterior except
that each of a finite number of points may belong to several contours.
\end{theo}

\begin{theo}
\label{PreobGaierTh1}
Let ${\mathcal C}$ be a closed limited point set whose complement $K$ is connected
and regular.
Suppose $\rho>1$ is the largest number such that $f$ is analytic inside $C_{\rho}$.
Choose $R_1$ and $R$ such that $1<R_1<R<\rho$.
Suppose $d_1$ and $d_2$ are such that
for $z\in{\mathcal C}$, $t_1\in C_{R_1}$ we have $d_1\le|t_1-z|\le d_2$.

Then there exists a set of points $\zeta_m^{(n)}$,
$n=1,2,\dots;$ $m=1,2,\dots,n$ \textup{(}the Fekete points of ${\mathcal C}$\textup{)} such that
for $z\in{\mathcal C}$, $t\in C_{R}$ we have an estimate
\[
\left|\frac{\omega_n(z)}{\omega_n(t)}\right|\le(n+2)\frac{d_2}{d_1}\left(\frac{R_1}{R}\right)^{n+1},\qquad
\omega_n(z)=(z-\zeta_1^{(n)})(z-\zeta_2^{(n)})\dots(z-\zeta_n^{(n)}).
\]
\end{theo}

\begin{theo}
\label{PreobWalshTh5}
Let ${\mathcal C}$ be a closed limited point set whose complement $K$ is connected
and regular. If the function $f(z)$ is single-valued and analytic on and within
$C_R$, there exists a sequence of polynomials $p_n(z)$ of respective degrees $n=0,1,2,\cdots$
such that we have
\[
|f(z)-p_n(z)|\le M(n+2)\frac{d_2}{d_1}\left(\frac{R_1}{R}\right)^{n+1},\qquad z\text{ on }{\mathcal C},
\]
where $M$ explicitly depends on $R$, but not on $n$ or $z$.
\end{theo}

\begin{figure}
  \begin{center}
    \includegraphics[width=\textwidth,keepaspectratio]{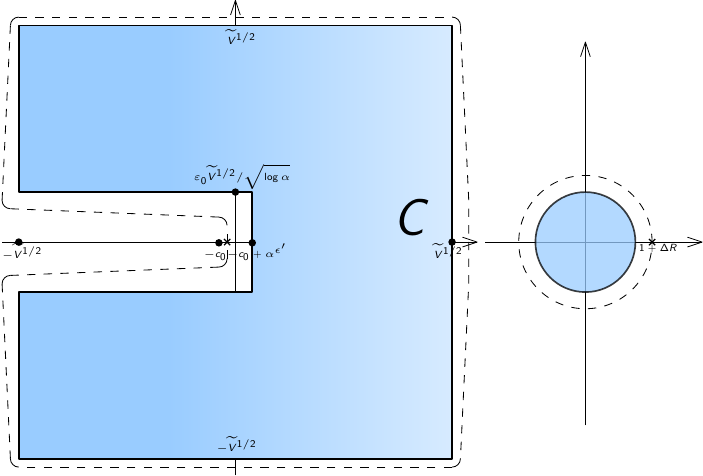}
    \caption{The set ${\mathcal C}$ and the image of $\phi(z)$} \label{PreobCFig}
  \end{center}
\end{figure}

To prove Theorem~\ref{PreobWalshTh5} with explicit constants, we need an explicit expression
for the function $w=\phi(z)$ defined in Theorem~\ref{PreobWalshTh1}
that maps the complement $K$ of our region~${\mathcal C}$
onto the exterior of the unit circle $\gamma$.
See Figure~\ref{PreobCFig}.
Such expression for the inverse function $z=\phi^{-1}(w)$
is supplied in a version of the Schwarz--Christoffel formula
(see~\cite[Chapter 6, \S\,2{.}2, Theorem~5]{PreobAhl79},
\cite[Section 47, Theorem~9{.}9]{PreobMar65}).

The main fact we need about the Schwarz--Christoffel function $\phi(z)$
is its asymptotic behavior in a neighborhood of $-c_0$.
This behavior is determined by the crowding phenomenon in Schwarz--Christoffel mapping
(see~\cite[Section 2{.}6]{PreobDT02}).

More precisely, let $\phi$ be defined in Theorem~\textup{\ref{PreobWalshTh1}}
and $D(-c_0,\alpha^{\eps'})$ be the disk of the $z$-plane centered at $-c_0$ of radius $\alpha^{\eps'}$.
Then for $\Delta R\asymp\frac{\alpha^{\eps'}}{\widetilde{V}^{1/2}}e^{-\sqrt{\log\alpha}/\eps_0}$
with $\eps'$ and $\eps_0$ going to $0$ sufficiently slowly as $\alpha$ goes to $\infty$,
the set
\[
\{z=\phi^{-1}(w)\ {:}\ |w|=1+\Delta R\}\cap D(-c_0,\alpha^{\eps'})
\]
lies within $\frac12\alpha^{\eps'}$ of the boundary $\partial{\mathcal C}$ in its exterior.

\textbf{Proof of Theorem~\ref{PreobWalshTh5}.}
Let $p_n(z)$ be the polynomial of degree $n$ which coincides
with $f(z)$ in the points
$\zeta_1^{(n+1)},\zeta_2^{(n+1)},\dots,\zeta_{n+1}^{(n+1)}$ of Theorem~\ref{PreobGaierTh1}.
Using~\cite[\S\,4{.}5, Equation~(11)]{PreobWal56} (Hermite's formula), we have
\[
f(z)-p_n(z)=\frac1{2\pi i}\inte_{C_R}\frac{\omega_{n+1}(z)f(t)\,dt}{\omega_{n+1}(t)(t-z)},\qquad z\in{\mathcal C}.
\]
Taking Theorem~\ref{PreobGaierTh1} into account, and recalling
\[
f(t)=\frac{\widetilde{\mathcal M}}{c_0(\alpha)+t},
\]
with $\widetilde{\mathcal M}\ll\alpha^{\eps'}$,
we see that for $z$ on ${\mathcal C}$ we have,
with $R_1=1+\frac12\Delta R$,
\[
|f(z)-p_n(z)|\ll\frac{n+2}{\alpha^{\eps'}\left(1+\frac12\Delta R\right)^n},\qquad z\text{ on }{\mathcal C}.
\]

The above construction needs to be adjusted in order to have the polynomials ${\mathcal L}_{n,{\mathcal C}}(z)$
with ${\mathcal L}_{n,{\mathcal C}}(0)=\frac1{c_0(\alpha)}$
which is necessary for Theorem~\ref{PreobLevCon} (this condition is essential in the application of Littlewood's lemma).

We choose the translation parameter $z_0=z_0(n)=O(\widetilde{V}^{1/2})$ and scale factor $a=a(n)=O(1)$
so that
\[
  {\mathcal L}_{n,{\mathcal C}}(z)=\frac1{\max_{z\ {:}\ az+z_0\in {\mathcal C}'\cup{\mathcal C}''}|p_n(az+z_0)|}
  p_n(az+z_0)
\]
satisfies ${\mathcal L}_{n,{\mathcal C}}(0)=\frac1{c_0(\alpha)}$.
The parameters exist by the Hadamard three-circle theorem.

Thus we have
\begin{equation}\label{PreobUniNorm}
  \begin{split}
    &\max_{z\ {:}\ az+z_0\in{\mathcal C}}\left(\left|{\mathcal L}_{n,{\mathcal C}}(z)\right|-\left|\frac{\widetilde{\mathcal M}}{c_0(\alpha)+az+z_0}\right|\right)
    \ll\frac{n+2}{\alpha^{\eps'}\left(1+\frac12\Delta R\right)^n},\\
    &\max_{z\ {:}\ az+z_0\in{\mathcal C}'\cup{\mathcal C}''}\left|{\mathcal L}_{n,{\mathcal C}}(z)\right|\le1.
  \end{split}
\end{equation}
with the absolute implied constant.

\section{Completion of proof of Theorem~\ref{Preob100percent}}\label{seccompletion}

Let in Theorem~\ref{PreobLevCon}
\[
a=\frac12-\frac{R-\delta_T}{\log T}
\]
with
\[
R=\epsilon\log\alpha,
\]
$\epsilon>0$ fixed,
and let
\[
\ds=\frac{\alpha}{\log T},
\]
where $\alpha>0$ is real and goes to infinity with $T$ sufficiently slowly.

Recall that in~\eqref{PreobUniNorm} $\widetilde{\mathcal M}$ is such that
\[
\widetilde{\mathcal M}\ll\alpha^{\epsilon'},
\]
and in the proof of Theorem~\ref{PreobExceptionTh} we established that we can take
\[
\widetilde{V}^{1/2}=\alpha^{C(1/2+2\epsilon)+c}
\]
with some fixed $c=o(C)>0$.

Let $\epsilon_0(\alpha)>0$ and $\epsilon'(\alpha)>0$ go to $0$ sufficiently slowly as $\alpha$ goes to $\infty$.

Using Theorems~\ref{PreobTranslationTh},~\ref{PreobApproximationTh},~\ref{PreobExceptionTh},
and denoting $E=E_{\alpha,\delta_T}\cup([T,2T]\setminus{\mathcal M}_\alpha)$
we can write
\[
\begin{split}
\frac2R\log I(R)&=\frac2R\log\left(\frac1{T-\eps_ET}\inte_{[T,2T]\setminus E}|F^{\ast}(a+it)|\,dt\right)\\
&=\frac2R\log\left(\frac1{T-\epsilon_ET}\int\limits_{[T,2T]\setminus E}|\zeta(a+\Delta\sigma+it)|
|c_0(\alpha)+\lambda(a+\Delta\sigma+it)|\right.\\
&\left.\vphantom{\int\limits_{[T,2T]\setminus E}\zeta}\times|M(a+it)||{\mathcal L}(A(a+\Delta\sigma+it))|\,dt\right),
\end{split}
\]
and
\[
\begin{split}
&\bigg|c_0(\alpha)+\lambda(a+\Delta\sigma+it)\bigg||{\mathcal L}(A(a+\Delta\sigma+it))|\\
&\ll\widetilde{\mathcal M}+\frac{\left|c_0(\alpha)+A(a+\Delta\sigma+it)\right|}{\widetilde{V}^{10}}
\end{split}
\]
by~\eqref{PreobSelbergUpperBound} and~\eqref{PreobUniNorm},
since we chose $n\asymp\frac{\widetilde{V}^{1/2}}{\alpha^{\eps'}}e^{\sqrt{\log\alpha}/\eps_0}\log\left(\widetilde{V}\right)$
and
\[
\max\limits_{t\in [T,2T]\setminus E}\left|\frac{\widetilde{\mathcal M}}{c_0+A(s)}-{\mathcal L}(A(s))\right|
\ll\frac1{\widetilde{V}^{10}}.
\]
We have $\frac{\log\widetilde{\mathcal M}}R$ arbitrarily small.
Now we apply the Cauchy--Schwarz inequality and show that for
\[
\bar{R}=\alpha-R+\delta(T)>0
\]
we have
\[
\inte_T^{2T}|\zeta(\frac12+\frac{\bar{R}}{\log T}+it)|^2
|M(a+it)|^2\,dt\le(1+\epsilon(\alpha,T))T
\]
with $\epsilon(\alpha,T)>0$ arbitrarily small.

To prove the estimate, in~\cite[Theorem 2]{PreobCon89} we take $R_C=-\bar{R}$
(note that our $\bar{R}>0$, but $R_C$ can be negative in Conrey's theorem),
$Q(x)=1$, and take $0<\theta<\frac47$ so that $\theta\bar{R}$ goes to infinity with $\alpha$.
Then Conrey's theorem asserts that for an optimal choice of $P$ in the mollifier $M(s)$ we have
(see~\cite[(49)]{PreobCon89})
\[
\inte_T^{2T}|\zeta(\frac12+\frac{\bar{R}}{\log T}+it)|^2
|M(a+it)|^2\,dt\thicksim c(1,\bar{R})T,
\]
with
\[
c(1,\bar{R})=\frac12\left(|w(0)|^2+|w(1)|^2\right)+A\bar{\alpha}\frac{\coth\bar{\alpha}}{\theta}.
\]
Here
\[
w(y)=\exp(-\bar{R}y),
\]
\[
A=\inte_0^1e^{-2\bar{R}y}dy,
\]
\[
\bar{\alpha}=\sqrt{\frac CA},
\]
\[
C=\theta^2\inte_0^1{\bar{R}}^2e^{-2\bar{R}y}dy,
\]
i.\,e.
\[
c(1,\bar{R})=\frac12+\frac12|w(1)|^2+\sqrt{A\frac C{\theta^2}}\coth\bar{\alpha}.
\]
We have
\[
w(1)=e^{-\bar{R}},
\]
\[
A=\frac1{2\bar{R}}\left(1-e^{-2\bar{R}}\right),
\]
\[
\frac C{\theta^2}=\frac{\bar{R}}2\left(1-e^{-2\bar{R}}\right),
\]
\[
A\frac C{\theta^2}=\frac14\left(1-e^{-2\bar{R}}\right)^2,
\]
\[
\bar{\alpha}=\theta\bar{R}.
\]
Thus, since $\bar{R}$ can be taken arbitrarily large, then $c(1,\bar{R})\le(1+\epsilon(\bar{R}))$,
with $\epsilon(\bar{R})$ arbitrarily small.

The remaining terms in Theorem~\ref{PreobLevCon} are small
by Theorem~\ref{PreobExceptionTh}.

This completes the proof of Theorem~\ref{Preob100percent}.

\medskip

\selectlanguage{english}

\bigskip

\end{document}